\documentclass[journal,final]{IEEEtran}

\usepackage[utf8]{inputenc}

\usepackage{etex}
\usepackage[vlined,ruled]{algorithm2e}
\usepackage[dvipsnames]{xcolor}
\usepackage{pgfpages}
\usepackage{amssymb, mathrsfs}
\usepackage{cite}

\usepackage[plainpages=false]{hyperref}
\hypersetup{
	colorlinks=true,       
	linkcolor=black,        
	citecolor=black,        
	breaklinks=true,
	urlcolor=blue,
}

\usepackage{dsfont}
\usepackage{siunitx}
\usepackage{verbatim}									
\usepackage{mathtools}									
\usepackage{siunitx}									
\usepackage{mdframed}

\usepackage{graphicx}
\usepackage{epstopdf}

\usepackage[font=small]{caption}
\usepackage{subcaption}
\usepackage[activate={true,nocompatibility},final,tracking=true,kerning=true,spacing=true,factor=1100,stretch=10,shrink=10]{microtype}

\usepackage{array}
\usepackage{multirow}
\usepackage{enumerate}
\usepackage{booktabs}

\graphicspath{{./fig/}}

\usepackage{pgfplots}
\pgfplotsset{width=8cm,compat=newest}
\newcommand*{\damping}{0.006}%
\newcommand*{\freq}{25}%
\pgfmathsetmacro{\freqd}{sqrt(1-(\damping)^2)*\freq}%

\pgfplotsset{
	standard/.style={
		axis x line=middle,
		axis y line=middle,
		enlarge x limits=0.15,
		enlarge y limits=0.15,
		every axis plot post/.style={mark options={fill=black}},
	}
}

\pgfplotsset{%
	,compat=1.12
	,every axis x label/.style={at={(current axis.right of origin)},anchor=north west}
	,every axis y label/.style={at={(current axis.above origin)},anchor=north east}
}

\usepackage{tikz}
\usetikzlibrary{arrows}
\usetikzlibrary{decorations.pathmorphing}
\usetikzlibrary{decorations.markings}
\usetikzlibrary{snakes}
\usetikzlibrary{matrix}
\usetikzlibrary{calc}
\usetikzlibrary{patterns}
\usetikzlibrary{positioning}
\usetikzlibrary{shapes}
\usetikzlibrary{shapes.symbols}
\usetikzlibrary{shapes.geometric}
\usetikzlibrary{shadows.blur}
\usetikzlibrary{shapes.arrows}
\usetikzlibrary{shapes.multipart}
\usetikzlibrary{shapes.callouts}
\usetikzlibrary{shapes.misc}

\tikzstyle{every node}=[font=\small]
\tikzstyle{every path}=[line width=0.8pt,line cap=round,line join=round]

\usepackage[american,smartlabels]{circuitikz}
\ctikzset{bipoles/thickness=1}
\ctikzset{bipoles/length=0.8cm}
\ctikzset{bipoles/diode/height=.375}
\ctikzset{bipoles/diode/width=.3}
\ctikzset{tripoles/thyristor/height=.8}
\ctikzset{tripoles/thyristor/width=1}
\ctikzset{bipoles/vsourceam/height/.initial=.7}
\ctikzset{bipoles/vsourceam/width/.initial=.7}

\newcommand{\real}{\mathbb{R}}

\newcommand{\setdef}[2]{\{#1 \;|\; #2\}}

\newcommand{\vect}[1]{\mathbbold{#1}}
\newcommand{\vones}[1][]{\vect{1}_{#1}}
\newcommand{\vzeros}[1][]{\vect{0}_{#1}}
\DeclareSymbolFont{bbold}{U}{bbold}{m}{n}
\DeclareSymbolFontAlphabet{\mathbbold}{bbold}

\newcommand{\map}[3]{#1: #2 \rightarrow #3}

\renewcommand{\theenumi}{(\roman{enumi})}

\newcommand{\tb}{\color{blue}}

\newcommand{\T}{\mathsf{T}} 

\newcommand\oprocendsymbol{\hbox{$\square$}}
\newcommand\oprocend{\relax\ifmmode\else\unskip\hfill\fi\oprocendsymbol}

\newtheorem{theorem}{Theorem}[section]
\newtheorem{lemma}[theorem]{Lemma}

\newtheorem{corollary}[theorem]{Corollary}

\newtheorem{remark}[theorem]{Remark}
\newtheorem{assumption}[theorem]{Assumption}
\newtheorem{example}[theorem]{Example}

\newenvironment{pfof}[1]{\vspace{1ex}\noindent{\itshape Proof of
		#1:}\hspace{0.5em}} {\hfill\oprocend\vspace{1ex}}
\newenvironment{proof}[1]{\vspace{1ex}\noindent{\itshape Proof:}\hspace{0.5em}} {\hfill\oprocend\vspace{1ex}}

\newcommand{\NI}{\mathsf{NI}}

\makeatletter
\newcommand{\vast}{\bBigg@{4}}
\newcommand{\Vast}{\bBigg@{5}}

\newcommand{\define}{\ensuremath{\triangleq}}

\newcommand{\ddiag}{\mathrm{diag}}
\newcommand{\bbm}{\begin{bmatrix}}
\newcommand{\ebm}{\end{bmatrix}}
\def\be{\begin{equation}}
\def\ee{\end{equation}}
\newcommand\ddfrac[2]{\frac{\displaystyle #1}{\displaystyle #2}}
\usepackage[normalem]{ulem}
\renewcommand{\tb}{\color{black}}
\begin{document}
\title{Diagonal Stability of Systems with Rank-1 Interconnections and Application to Automatic Generation Control in Power Systems}
\author{John W. Simpson-Porco and Nima Monshizadeh 
\thanks{This work was supported in part by NSERC Discovery Grant RGPIN-2017-04008.}
\thanks{J.~W.~Simpson-Porco is with the Department of Electrical and Computer Engineering, University of Toronto, 10 King's College Road,
Toronto, ON, M5S 3G4, Canada (email: jwsimpson@ece.utoronto.ca).}
\thanks{N.~Monshizadeh is with the Engineering and Technology
Institute, University of Groningen, 9747AG, The Netherlands (e-mail:
n.monshizadeh@rug.nl).}}

\maketitle

\begin{abstract}
We study a class of matrices with a rank-1 interconnection structure, and derive a simple necessary and sufficient condition for diagonal stability. The underlying Lyapunov function is used to provide sufficient conditions for diagonal stability of approximately rank-1 interconnections. The main result is then leveraged as a key step in a larger stability analysis problem arising in power systems control. Specifically, we provide the first theoretical stability analysis of automatic generation control (AGC) in an interconnected nonlinear power system. {\tb Our analysis is based on singular perturbation theory, and provides theoretical justification for the conventional wisdom that AGC is stabilizing under the typical time-scales of operation. We illustrate how our main analysis results can be leveraged to provide further insight into the tuning and dynamic performance of AGC.}
\end{abstract}


%


\section{Introduction}
\label{Sec:Introduction}

Classical Lyapunov theory for linear systems states that eigenvalues of a matrix $A$ are contained in the open left-half complex plane if and only if there exists a positive definite matrix $D \succ 0$ such that $A^{\T}D + DA \prec 0$; such a matrix $A$ is called (Hurwitz) stable \cite[Thm. 2.2.1]{RAH-CRJ:94}. In general, solutions $D$ to this inequality will be dense matrices, where most or all matrix entries are non-zero. {\tb In many applications, listed in the sequel,}  it is desirable to further find a solution $D \succ 0$ which is a diagonal matrix. If this can be done, then the matrix $A$ is \emph{diagonally stable} \cite{EK-AB:00}.

{\tb While any diagonally stable matrix is stable, the converse is false; diagonal stability is indeed a considerably stronger property than Hurwitz stability, due to the removal of degrees of freedom from the Lyapunov matrix $D$.} With this strength however comes significant advantages, and diagonal stability has found widespread application in areas such as economics \cite{EK-AB:00}, biological systems \cite{LS-MA-EDS:10},  singular perturbation theory \cite{HKK:87}, positive systems analysis \cite{RS-KSN:09,AR:15}, and analysis of general large-scale interconnected systems \cite{MH:78,DDS:78,MV:81,MA-CM-AP:16}. A significant body of literature exists delineating classes of diagonally stable matrices. Prominent examples of such classes include symmetric negative definite matrices, Hurwitz lower/upper triangular matrices, $M$-matrices, and certain types of cyclic interconnections \cite{MArcak:11}; see \cite{EK-AB:00} for even more nuanced classes of matrices.

 Our initial focus in this paper is to further contribute to the theory of diagonal stability by presenting necessary and sufficient conditions for diagonal stability of a new class of matrices, consisting of rank-1 perturbations of negative definite diagonal matrices. Our motivation stems from the fact that this class of matrices arises in stability analysis of certain networked systems, such as automatic generation control (AGC) in interconnected power systems, and diagonal stability of such a matrix is precisely the condition required to complete a Lyapunov-based stability analysis. The second half of this paper consists of a detailed and self-contained treatment of the AGC application with its corresponding introduction and literature deferred to Section \ref{Sec:BackgroundAGC}. 
 
 As an independent system-theoretic motivation for the study of such a class of matrices, we begin with an example arising in the stability analysis of interconnected nonlinear systems. Consider a collection of $N \geq 2$ single-input single-output nonlinear systems
 \begin{equation}\label{Eq:NonlinearSISO}
\begin{aligned}
\dot{x}_i &= f_i(x_i,u_i)\\
y_i &= h_i(x_i,u_i),
\end{aligned} \qquad i \in \{1,\ldots,N\},
\end{equation}
with internal states $x_i \in \real^{n_i}$, inputs $u_i \in \real$ and outputs $y_i \in \real$. The functions $f_i$ and $h_i$ are sufficiently smooth on a domain containing the origin, and satisfy $f_i(0,0) = 0$ and $h_i(0,0) = 0$ for all $i \in \{1,\ldots,N\}$. Assume that each subsystem is \emph{output-strictly passive} (see e.g., \cite{MA-CM-AP:16}), meaning that there exists a constant $\delta_i > 0$ and a continuously differentiable storage function $\map{V_i}{\real^{n_i}}{\real}$ which is positive definite with respect to the origin $x_i = 0$ and satisfies the dissipation inequality
\begin{equation}\label{Eq:OSP}
\nabla V_i(x_i)^{\T}f_i(x_i,u_i) \leq -\delta_i y_i^2 + y_i u_i
\end{equation}
on some neighbourhood of $(x_i,u_i) = (0,0)$. Suppose now that the subsystems in \eqref{Eq:NonlinearSISO} are interconnected according to 
\begin{equation}\label{Eq:FeedbackLaw}
u_i = k_i \sum_{j=1}^{n}\nolimits g_j y_j, \qquad i \in \{1,\ldots,N\},
\end{equation}
where $k_i, g_j \in \real$ are constants. {We interpret \eqref{Eq:FeedbackLaw} as a ``gather-and-broadcast'' type control law: the local measurements $y_j$ are centrally collected, the linear combination $\eta \define \sum_{j=1}^{N}g_jy_j$ is centrally computed using those local measurements and broadcast to each subsystem, and the local inputs are then determined as $u_i = k_i \eta$. {\tb Control architectures of this form arise in optimal frequency control \cite{FD-SG:17}, \cite{NM-CDP-AJvdS-JMAS:18}, Nash seeking algorithms in aggregative games \cite{CDP-SG:20},  and in network flow control \cite{XF-TA-MA-JTW-TB:06,JTW-MA:04}.} A simple example would be the case where $g_j = \tfrac{1}{N}$ and $k_i = 1$ in which each subsystem is driven by the arithmetic average of all subsystem output signals.}
 
To assess stability of the origin of the interconnected system \eqref{Eq:NonlinearSISO},\eqref{Eq:FeedbackLaw}, one may follow the well-established (e.g., \cite{MA-CM-AP:16}) approach of defining a composite Lyapunov candidate
\[
V(x_1,\ldots,x_N) = \sum_{i=1}^{N}\nolimits d_iV_i(x_i)
\]
with coefficients $d_i > 0$ to be determined. Using the dissipation inequality \eqref{Eq:OSP}, the derivative of $V$ along \eqref{Eq:NonlinearSISO},\eqref{Eq:FeedbackLaw} can now be computed as
{\tb
\[
\begin{aligned}
&\dot{V} = \sum_{i=1}^{N}\nolimits d_i \nabla V_i(x_i)^{\T}f_i(x_i, u_i)\\ & \; \leq  \sum_{i=1}^{N}\nolimits d_i(-\delta_i y_i^2 + y_iu_i)= \sum_{i=1}^{N} d_i(-\delta_i y_i^2 + y_i k_i \sum_{j=1}^{N} g_jy_j),
\end{aligned}
\]
}
which in vector notation can be compactly written as
{\tb
\[
\dot{V} \leq \boldsymbol{y}^{\T}DA\boldsymbol{y} = \tfrac{1}{2}\boldsymbol{y}^{\T}(DA+A^{\T}D)\boldsymbol{y}
\]
}
where $\boldsymbol{y} = \mathrm{col}(y_1,\ldots,y_N)$, $D = \mathrm{diag}(d_1,\ldots,d_N)$ and 
\begin{equation}\label{Eq:AStabilityAnalysis}
A = -\mathrm{diag}(\delta_1,\ldots,\delta_N) +  \boldsymbol{k}\boldsymbol{g}^{\T}
\end{equation}
with 
$\boldsymbol{k}=\mathrm{col}(k_1,\ldots,k_N)$ and $\boldsymbol{g} = \mathrm{col}(g_1,\ldots,g_N)$.\footnote{See Section \ref{ss:notation} for the notation used.} Note that $A$ is a rank-1 perturbation of a negative definite diagonal matrix. If $A$ is diagonally stable, then there will exist $\epsilon > 0$ and a selection of $(d_1,\ldots,d_N)$ such that $A^{\T}D+DA \preceq -2\epsilon I_N$. This will in turn imply that $\dot{V} \leq -\epsilon \|\boldsymbol{y}\|_2^2$, and local asymptotic stability of the origin can then be proven under standard zero-state detectability conditions on each subsystem in \eqref{Eq:OSP} \cite{HKK:02}. The key step in the analysis is therefore to develop conditions under which \eqref{Eq:AStabilityAnalysis} is diagonally stable.

\smallskip

\subsection{Statement of Contributions} This paper contains two main contributions. First, in Section \ref{Sec:Main} we state and prove necessary and sufficient conditions for diagonal stability of a matrix consisting of a negative definite diagonal matrix plus a rank-1 perturbation under certain sign constraints (Theorem \ref{Thm:Main} and Corollary \ref{Cor:Explicit}). To the best of our knowledge,  diagonal stability of this class of matrices has not been previously studied. We illustrate how our conditions are stronger than those required for Hurwitz stability, but weaker than those imposed by a small-gain approach. Finally, we extend our analysis to provide a sufficient condition for diagonal stability of approximately rank-1 interconnections. The results are illustrated via a simple numerical example.

Second, in Section \ref{Sec:AGC}, we leverage our results from Section \ref{Sec:Main} as a key step in a larger detailed stability treatment of AGC in interconnected power systems. Despite being one of the first large-scale control systems deployed in industry, to the best of our knowledge no general dynamic stability analysis of AGC is available in the literature. The standard textbook analyses \cite{AJW-BFW:96,PK:94} are based on equilibrium analysis, and do not consider dynamic stability of the equilibrium. Notable exceptions to a purely static analysis are the treatments in \cite{MDI-SL:96, MI-JZ:00}, but these analyses focus on local reduced-order area models, and do not analyze the interconnected dynamics of power systems with AGC. We provide here the first nonlinear stability treatment of AGC, keeping the treatment reasonably self-contained, and accessible to those with only passing familiarity with power systems control. {Readers interested only in this application can skip directly to Section \ref{Sec:AGC}, with minimal loss of continuity.} Roughly speaking, our main result (Theorem \ref{Thm:AGCStable}) states that under usual time-scales of operation, power systems with AGC are internally stable {\tb for any tuning of the frequency bias gains used within the AGC design. This result provides a rigorous theoretical backing for the practical observation that AGC is stable in practice, and for the accepted engineering practice of decentralized and uncoordinated tuning of AGC systems. Finally, in Section \ref{Sec:DynamicPerformance} we examine some of the implications of our main analysis result for the dynamic performance of AGC systems. In particular, we demonstrate that the industrial tuning practice for AGC \cite{NERC:10bal} can be rigorously understood within our analysis framework, and that our framework has the potential for generating other novel insights into the dynamic performance of AGC systems.}

\smallskip

\subsection{Notation}\label{ss:notation}
The identity matrix of size $n$ is denoted by $I_n$, and $\vones[n]$ denotes the vector of all ones in $\real^n$.  We denote the set of vectors in $\real^n$ with all positive, nonnegative, and nonzero elements, by $\real^n_{>0}$, $\real^n_{\geq 0}$, $\real^n_{\neq 0}$, respectively. The set of numbers $\{1, \ldots, N\}$ is denoted by $\mathcal{I}_N$. Given an ordered set of scalars or column vectors $(x_1,\ldots,x_n)$, we let $\mathrm{col}(x_1,\ldots,x_n)$ denote the concatenated column vector. As usual, $A \succ 0$ means $A$ is positive definite. We use $\ddiag(x_i)$ to denote the diagonal matrix with its $i$th diagonal element equal to  
$x_i$. We sometimes denote the latter by $\ddiag(x)$ as well for a given vector $x$. In the same vein, we use $\ddiag(x,y)$ to denote the diagonal matrix
\[
\bbm
\ddiag (x) & 0\\
0 & \ddiag(y)
\ebm
\]
for two vectors $x$ and $y$.
For a scalar $a$, we define the operator 
\[
[a]_{+}:=\max\, (a, 0).
\]
which truncates negative values to zero. We use the same notation $[a]_{+}$ for a vector $a$ by interpreting the $\max$ operator element-wise with respect to a zero vector.  



\section{Diagonal Stability of a Diagonal Matrix with a Rank-1 Perturbation}
\label{Sec:Main}

We begin by considering $N \times N$ matrices of the form
\begin{equation}\label{Eq:A-s}
A=- \Delta + S 
\end{equation}
where $\Delta =\mathrm{diag}(\delta_1,\ldots,\delta_N) \succ 0$ and $S \in \real^{N\times N}$. In an interconnected systems framework, the matrix $\Delta$ captures the dissipation of each local subsystem, and the matrix $S$ represents the interconnection structure among the subsystems. The matrix \eqref{Eq:A-s} can appear directly as the state matrix of a linear system, or arise as part of a quadratic form in a composite Lyapunov stability analysis of a nonlinear interconnected system as discussed in Section \ref{Sec:Introduction}.
We are interested in determining conditions under which \eqref{Eq:A-s} is diagonally stable.

In the absence of additional structural information on $S$, the standard method to enforce diagonal stability of $A$ is a {\em small gain} approach. This amounts to assuming that the norm of $S$ is sufficiently small compared to $\Delta$, yielding the condition $\|\Delta^{-1}S\|_2 < 1$, in which case $A$ is diagonally stable with $D = I_N$. When additional structural information is available, less conservative conditions can be obtained. For example, if $S$ is a non-negative matrix, then diagonal stability of $A$ can be studied using the theory of M-matrices, resulting in non-conservative conditions \cite{EK-AB:00}. Another notable example of a tight diagonal stability condition is the secant condition for systems with a cyclic interconnection structure \cite{MArcak:11}.     

Our primary case of interest is that where $\mathrm{rank}(S)=1$, in which case the matrix $S$ can be written as
\begin{equation}\label{Eq:S}
S= x y^\T
\end{equation}
for some vectors $x, y \in \real^N$.
%
%
As diagonal stability is preserved under left or right multiplication by any positive definite diagonal matrix {\tb \cite[Lem. 2.1.4]{EK-AB:00}}, it follows that diagonal stability of $A = -\Delta + x  y^{\T}$ is equivalent to diagonal stability of 
\begin{equation}\label{Eq:A}
A_I \define -I_N + \hat x y^{\T},
\end{equation}
where $\hat x:=\Delta^{-1} x$. For matrices of the form \eqref{Eq:A}, the first question is whether there is a gap between Hurwitz stability and diagonal stability. Note that the $A_I$ has $N-1$ eigenvalues at $-1$, with its remaining eigenvalue at $-1 + \hat x^{\T}y$. It follows that $A_I$ is Hurwitz stable if and only if
\begin{equation}\label{Eq:HurwitzCondition}
\sum_{i=1}^{N}\nolimits \hat{x}_iy_i < 1.
\end{equation}
The inequality in \eqref{Eq:HurwitzCondition}, however, does not ensure diagonal stability of \eqref{Eq:A} {\tb as illustrated in the following example.}

\smallskip

\begin{example}\label{ex:num}
Consider $A_I \in \real^{2\times 2}$ given by
\begin{equation}\label{Eq:Example}
A_I = -I_2 + \alpha \begin{bmatrix}1 \\ -1\end{bmatrix}\begin{bmatrix}1\\ 1\end{bmatrix}^{\T} = \begin{bmatrix}
-1+\alpha & \alpha\\
-\alpha & -1-\alpha
\end{bmatrix},
\end{equation}
for $\alpha \in \real$.
The Hurwitz stability condition \eqref{Eq:HurwitzCondition} is satisfied independent of $\alpha$, since in this case {\tb $\sum_{i=1}^{N}\hat{x}_iy_i = \alpha\left[\begin{smallmatrix}1 & -1\end{smallmatrix}\right]\left[\begin{smallmatrix}1 \\ 1\end{smallmatrix}\right] = 0$ for all $\alpha \in \real$}. Diagonal stability imposes the existence of $D = \mathrm{diag}(d_1,d_2) \succ 0$ such that
\begin{equation}\label{Eq:AIExample}
A_I^{\T}D + D A_I = \begin{bmatrix}
2d_1(-1+\alpha)  & \alpha(d_1-d_2)\\
\alpha(d_1-d_2) & -2d_2 (1+\alpha)
\end{bmatrix} \prec 0.
\end{equation}
A necessary condition for this is $|\alpha|<1$, in which case \eqref{Eq:AIExample} holds if and only if
\[
\alpha^2\left(\tfrac{d_1}{d_2}\right)^2- (4-2\alpha^2)  \left(\tfrac{d_1}{d_2}\right)+\alpha^2<0.
\]
There exist $d_1, d_2>0$ satisfying above if and only if $(4-2\alpha^2)^2-4\alpha^4>0$, which holds if and only if $|\alpha|<1$. Consequently \eqref{Eq:Example} is diagonally stable if and only if $|\alpha|<1$, whereas it is Hurwitz for all $\alpha\in \real$.
\oprocend
\end{example}

\medskip

In general, diagonal stability of \eqref{Eq:A} could be enforced through the previously mentioned small-gain criteria, yielding the condition
\begin{equation}\label{Eq:SmallGain}
\|\hat{x}y^{\T}\|_2 = \|\hat{x}\|_2 \|y\|_2 < 1.
\end{equation}
The condition \eqref{Eq:SmallGain} however is unnecessarily conservative, as it neglects the sign patterns (i.e., phase) of $\hat{x}$ and $y$. For Example \ref{ex:num}, \eqref{Eq:SmallGain} yields $|\alpha| < \tfrac{1}{2}$ instead of the necessary and sufficient condition $|\alpha| < 1$. A more dramatic example is the case $y=-\hat x$, where the resulting matrix $A_I  = -I_N - \hat x \hat x^\T$ is diagonally stable, independent of the norm of the vector $\hat{x}$. 


We now come to our main result, which avoids the small-gain type conditions as in \eqref{Eq:SmallGain}, and strengthens the Hurwitz stability condition \eqref{Eq:HurwitzCondition} to obtain a necessary and sufficient condition for diagonal stability. We require the assumption that either $x$ or $y$ is a nonnegative vector; this yields additional structure in the sign pattern of $xy^{\T}$. Since diagonal stability of $A$ and $A^\T$ are equivalent, we will simply assume that $y\in \real^N_{\geq 0}$. 

\medskip{}

\begin{theorem}\label{Thm:Main}
Let $A$ be given by \eqref{Eq:A-s} and \eqref{Eq:S}, with $x \in \real^{N}$, $y \in \real^{N}_{\geq 0}$.
Then $A$ is diagonally stable if and only if
\begin{equation}\label{Eq:DiagStabCondition}
\sum_{i=1}^{N}\nolimits \frac{1}{\delta_i}[x_iy_i]_{+} < 1.
\end{equation}
\end{theorem}


\smallskip

\begin{proof}
\emph{Sufficiency:} 
By defining $\hat x:=\Delta^{-1}x$, we equivalently show diagonal stability of $A_I$ in \eqref{Eq:A}.
If $x = \vzeros[N]$ or $y = \vzeros[N]$ then the conclusion holds, so assume that $x,y \neq \vzeros[N]$, and without loss of generality assume that the elements of $\hat x$ and $y$ are ordered such that $\hat x = \mathrm{col}(\hat x_1, \hat x_2)$ and $y = \mathrm{col}(y_1,\vzeros[N-r])$ with $y_1 \in \real_{>0}^{r}$, $\hat x_1 \in \real^{r}$, $\hat x_2 \in \real^{N-r}$, and $r \in \{1,\ldots,N\}$. With this re-ordering, we have that
\[
A_I = -I_{N} + \left[\begin{smallmatrix}\hat x_1 \\ \hat x_2\end{smallmatrix}\right]\left[\begin{smallmatrix}y_1 \\ \vzeros[N-r]\end{smallmatrix}\right]^{\T} = \begin{bmatrix}-I_{r} + \hat{x}_1y_1^{\T} & \vzeros[]\\
\hat{x}_2y_1^{\T} & -I_{N-r}
\end{bmatrix}.
\]
{\tb 
Note that $A_I$ is diagonally stable if and only if $A_I^{\T}$ is diagonally stable \cite[Lem. 2.1.7] {EK-AB:00}. Applying Lemma \ref{Lem:1} to $A_I^{\T}$, we find that} $A_I$ is diagonally stable if and only if $A_1 \define -I_r + \hat x_1y_1^{\T}$ is diagonally stable. We now apply the previous argument again to $\hat x_1 \in \real^r$. If $\hat{x}_1 = \vzeros[r]$ then the conclusion holds, so assume $\hat{x}_1 \neq \vzeros[r]$, and without loss of generality re-order the elements of $\hat{x}_1$ and $y_1$ such that $\hat x_1 = \mathrm{col}(\hat{x}_{+},-\hat{x}_{-},\vzeros[r-p-q])$ and $y_1 = \mathrm{col}(y_{+},y_{-},y_{\rm z})$ for some strictly positive vectors $\hat{x}_{+},y_{+} \in \real^p$, $\hat{x}_{-},y_{-} \in \real^{q}$ and $y_{\rm z} \in \real^{r-p-q}$ where $p+q \in \{1,\ldots,r\}$.
%
%
%
In this notation, we may write
\begin{equation}\label{Eq:AExpanded}
A_1 = \left[\begin{array}{c|c}
-I_{p+q} + \left[\begin{smallmatrix}\hat{x}_{+}\\-\hat{x}_{-}\end{smallmatrix}\right]\left[\begin{smallmatrix}y_{+}\\y_{-}\end{smallmatrix}\right]^{\T} & \left[\begin{smallmatrix}\hat{x}_{+}\\-\hat{x}_{-}\end{smallmatrix}\right]\left[\begin{smallmatrix}y_{\rm z}\end{smallmatrix}\right]^{\T}\\
\hline
\vzeros[] & -I_{r-p-q}
\end{array}\right].
\end{equation}
By Lemma \ref{Lem:1} again, $A_1$ is diagonally stable if and only if
\be\label{Eq:A11-1}
A_{11} \define -I_{p+q} + \left[\begin{smallmatrix}\hat{x}_{+}\\-\hat{x}_{-}\end{smallmatrix}\right]\left[\begin{smallmatrix}y_{+}\\y_{-}\end{smallmatrix}\right]^{\T}
\ee
is diagonally stable. We rewrite ${A}_{11}$ as 
\be\label{Eq:A11}
A_{11}=
\ddiag(\hat{x}_+,\hat{x}_-)
\underbrace{\big(-M+ \left[\begin{smallmatrix}\vones[p] \\ -\vones[q]\end{smallmatrix}\right]\left[\begin{smallmatrix}\vones[p] \\ \vones[q]\end{smallmatrix}\right]^{\T}
\big)}_{:=\overline A_{11}}
\ddiag(y_+,y_-)
\ee
where $M \in \real^{(p+q)\times(p+q)}$ is given by
 \begin{equation}\label{Eq:P}
 M = \left[\begin{smallmatrix}M_{11} & \vzeros[]\\
\vzeros[] & M_{22}\end{smallmatrix}\right] = \mathrm{diag}(\hat{x}_+,\hat{x}_-)^{-1}\mathrm{diag}(y_+,y_-)^{-1},
 \end{equation}
Clearly $M$ is diagonal and $M \succ 0$. 
By Lemma \ref{Lem:properties}, the matrix $A_{11}$ is diagonally stable if and only if $\overline A_{11}$ is diagonally stable. A direct calculation now shows that
\be\label{Eq:lyap}
\overline {A}_{11}^{\T}+\overline {A}_{11} = -2\overline{Q}_{11}
\ee
where
\begin{equation}\label{Eq:Q11bar}
\overline{Q}_{11} \define \begin{bmatrix}M_{11} - \vones[p]\vones[p]^{\T} & \vzeros[] \\ \vzeros[] & M_{22} + \vones[q]\vones[q]^{\T}\end{bmatrix}.
\end{equation}
If $p = 0$, then $\overline{Q}_{11} \succ 0$, and the conclusion follows. If $p \geq 1$, then observe that $\overline{Q}_{11} \succ 0$ if and only if $M_{11} - \vones[p]\vones[p]^{\T} \succ 0$. Performing a congruence transformation with $M_{11}^{-1/2} \succ 0$, we have that
\begin{equation}\label{Eq:SuffFinalInequality}
M_{11} - \vones[p]\vones[p]^{\T} \succ 0\,\,\, \Longleftrightarrow \,\,\, I_{p} - M_{11}^{-1/2}\vones[p]\vones[p]^{\T}M_{11}^{-1/2} \succ 0.
\end{equation}
The second term in the latter inequality is a rank-1 matrix, with unique positive eigenvalue given by
\[
\begin{aligned}
\mathrm{Trace}(M_{11}^{-1}) &= (\hat{x}_{+})^{\T}y_+ = \sum_{i=1}^{p}\nolimits \hat x_i y_i = \sum_{i=1}^{N}\nolimits [\hat x_iy_i]_{+}\\
&= \sum_{i=1}^{N}\nolimits \tfrac{1}{\delta_i}[x_iy_i]_{+} < 1
\end{aligned}
\]
due to \eqref{Eq:DiagStabCondition}. it follows that \eqref{Eq:SuffFinalInequality} holds, and therefore $\overline{A}_{11}$ is diagonally stable with certificate $\overline{D}_{11} = I_{p+q}$, which shows the result.

\medskip{}
\emph{Necessity:} 
Suppose that $A$ (and hence, $A_I$) is diagonally stable. Then any principal submatrix of $A_I$ is diagonally stable \cite[Lem. 2.1.8]{EK-AB:00}, and in particular then, the matrix $A_{11}$ in \eqref{Eq:A11} is diagonally stable. It follows next from Lemma \ref{Lem:properties} that $\overline A_{11}$ is also diagonally stable. Therefore,
the principal submatrix 
\[
-\ddiag(\hat{x}_+)^{-1}\ddiag(y_+)^{-1}+\vones[p]\vones[p]^\T,
\]
of $\overline A_{11}$ is diagonally stable as well. Consequently, by Lemma \ref{Lem:properties}, the matrix $-I_p+\hat{x}_+ y_+^\T$ is diagonally stable and thus Hurwitz. From \eqref{Eq:HurwitzCondition}, Hurwitz stability implies that
\[
(\hat{x}_{+})^{\T}y_+ = \sum_{i=1}^{p}\nolimits \hat{x}_{+,i}y_{+,i} = \sum_{i=1}^{N}\nolimits \frac{1}{\delta_i} [x_i y_i]_{+} < 1,
\]
which is precisely \eqref{Eq:DiagStabCondition}, completing the proof.
\end{proof}

\bigskip
It is illustrative to examine the result of Theorem \ref{Thm:Main} when $\Delta=I_N$, in which case we can compare the necessary and sufficient conditions for Hurwitz stability \eqref{Eq:HurwitzCondition} and diagonal stability \eqref{Eq:DiagStabCondition}. In contrast with \eqref{Eq:HurwitzCondition}, the condition \eqref{Eq:DiagStabCondition} for diagonal stability places tighter restrictions on the elements of $x$ and $y$. In particular, \eqref{Eq:HurwitzCondition} permits negative elements to compensate positive elements in the sum, whereas \eqref{Eq:DiagStabCondition} drops all negative elements from consideration. Revisiting Example \ref{ex:num}, we see that in this example the condition \eqref{Eq:DiagStabCondition} readily gives
\[
[\alpha]_++[-\alpha]_+=|\alpha|<1,
\]
which coincides with condition derived in the example.

If it happens that neither $x$ nor $y$ contain zero elements, then the diagonal stability certificate in Theorem \ref{Thm:Main} admits the  explicit form below.\footnote{Again note that an analogous result holds for $y\in \real^N_{\neq 0}$ and $x\in \real^N_{>0}$.}

\medskip

\begin{corollary}\label{Cor:Explicit}
Let $A$ be given by \eqref{Eq:A-s} and \eqref{Eq:S}, with $x \in \real^N_{\neq 0}$, $y \in \real^{N}_{> 0}$. Suppose that \eqref{Eq:DiagStabCondition} holds. Then with
\[
D=\ddiag(d_1,\ldots,d_N), \quad d_i= \frac{y_i}{|x_i|}, \quad i \in \{1,\ldots,N\}
\]
it holds that $A^{\T}D + DA \preceq -2\mu \Delta D$, 
where
\be\label{Eq:alpha}
\mu:= 1-\sum_{i=1}^{N}\nolimits \frac{1}{\delta_i}[x_iy_i]_{+} > 0.
\ee
\end{corollary}


\medskip{}
\begin{proof}
With $A = -\Delta + xy^{\T}$, we have that $A_{I} \define \Delta^{-1}A = -I_N + \hat{x}y^{\T}$ with $\hat{x} = \Delta^{-1}x$. Since all elements of $\hat{x}$ are non-zero and all elements of $y$ are positive, there exists a permutation matrix $\Pi$ such that
\[
\Pi^{\T} {\hat x} = \mathrm{col}(\hat{x}_{+},-\hat{x}_{-}), \quad \Pi^{\T} y = \mathrm{col}(y_{+},y_{-})
\]
where all subvectors $\hat{x}_+$, $\hat{x}_-$, $y_+$, $y_-$ are strictly positive. It follows that
\[
\Pi^{\T}A_{I}\Pi = A_{11} \define -I_{N} + \left[\begin{smallmatrix}\hat{x}_{+}\\-\hat{x}_{-}\end{smallmatrix}\right]\left[\begin{smallmatrix}y_{+}\\y_{-}\end{smallmatrix}\right]^{\T}.
\]
The matrix $A_{11}$ can be written as in \eqref{Eq:A11}
where
\[
\overline{A}_{11} = -\mathrm{diag}(\hat{x}_+,\hat{x}_-)^{-1}\mathrm{diag}(y_+,y_-)^{-1} + \left[\begin{smallmatrix}\vones[p] \\ -\vones[q]\end{smallmatrix}\right]\left[\begin{smallmatrix}\vones[p] \\ \vones[q]\end{smallmatrix}\right]^{\T}.
\]
We have that the matrix $\overline{A}_{11}$ satisfies \eqref{Eq:lyap} with $\overline{Q}_{11}$ given by \eqref{Eq:Q11bar}. 
Working backwards now, by Lemma \ref{Lem:properties}, the stability certificate for $A_{11}$ is given by
\be\nonumber
D_{11}:= \ddiag(\hat{x}_{+},\hat{x}_{-})^{-1}\ddiag(y_{+},y_{-}),
\ee
and in particular we have
\[
A_{11}^{\T}D_{11} + D_{11}A_{11} = -2\underbrace{\ddiag(y_{+},y_{-})\overline{Q}_{11}\ddiag(y_{+},y_{-})}_{\define Q_{11}}.
\]
To obtain a simple bound on $Q_{11}$, one may use \eqref{Eq:P} and repeat again the argument following \eqref{Eq:Q11bar} to show that
\begin{align*}
Q_{11}&= \ddiag(y_{+},y_{-}) \left(M -
\left[
\begin{smallmatrix}
\vones[p]\vones[p]^\T & 0\\
0 & -\vones[q]\vones[q]^\T
\end{smallmatrix}\right]\right)
\ddiag(y_{+},y_{-})\\
&\succeq \mu \,
\ddiag(y_{+},y_{-})\,  M \, \ddiag(y_{+},y_{-})\\
&= \mu D_{11}
\end{align*}
where $\mu$ is given by \eqref{Eq:alpha}. Since $A_I = \Pi A_{11} \Pi^{\T}$, we may again apply Lemma \ref{Lem:properties}, with
\[
D_{I} := \Pi D_{11}\Pi^{\T} = \ddiag(y_i/|\hat{x}_i|) = \ddiag(\delta_i y_i/|x_i|)
\]
to obtain
\[
\begin{aligned}
A_{I}^{\T}D_{I} + D_{I}A_{I} &= -2\Pi Q_{11}\Pi^{\T}.
\end{aligned}
\]
Finally, with $D = D_{I}\Delta^{-1}  = \ddiag(y_i/|x_i|)$ we have that
\[
\begin{aligned}
A^{\T}D + DA &= A_I^{\T}D_I + D_IA_I = -2 \Pi Q_{11} \Pi^{\T} \\ 
&\quad \preceq -2\mu  \Pi D_{11} \Pi^\T
= -2 \mu D_I
= -2 \mu \Delta D, 
\end{aligned}
\]
which completes the proof.
\end{proof}

\smallskip

The explicit certificate in Corollary \ref{Cor:Explicit} can be further leveraged to provide a sufficient condition for diagonal stability when the interconnection matrix in \eqref{Eq:A-s} is a perturbation of a rank one matrix. Indeed, suppose that $S$ now has the form
\[
S = xy^{\T} + \sigma E
\]
where $x \in \real^n_{\neq 0}$, $y \in \real^{n}_{>0}$,  $\sigma \in \real$, and $E \in \real^{n \times n}$; without loss of generality, we may assume that $\|E\|_2 = 1$. Assume that \eqref{Eq:DiagStabCondition} holds, and let $A_1 = -\Delta + xy^{\T}$, in which case we may write $A \define A_1 + \sigma E$. The conditions of Corollary \ref{Cor:Explicit} are satisfied, and we use the certificate $D \succ 0$ to compute that
\begin{align*}
A^\T D+D A&=A_1^\T D+D A_1 + \sigma (E^\T D + DE)\\ 
&\preceq -2\mu \Delta D +\sigma (E^\T D + DE)\\
&\preceq -2\mu \Delta D +2|\sigma|\|D\|_2 I_N\\
&\preceq -2\left(\mu (\min_{i}\delta_id_i) - |\sigma|(\max_{i}d_i)\right)I_N,
\end{align*}
from which we conclude that $A$ is diagonally stable with certificate $D$ if
\begin{equation}\label{Eq:SigmaBound}
|\sigma| < \mu \cdot \frac{\min_{i\in\mathcal{I}_N}d_i\delta_i}{\max_{i\in\mathcal{I}_N}d_i}, \quad d_i=\frac{y_i}{|x_i|}.
\end{equation}
Inequality \eqref{Eq:SigmaBound} is a robustness result, which restricts the norm of the perturbation. We illustrate the ideas with a simple example. 

\medskip

\begin{example}
Let $a \in \real$ and consider the matrix
\[
\begin{aligned}
A &= \begin{bmatrix}2a-8 & -a-3 & -3 & -a-3 & -a-3\\ 1-a & -3 & 1-a & a+1 & 1-a\\ -1 & a-1 & -a-4 & -1 & -1\\ -3 & a-3 & -a-3 & a-7 & -3\\ a+2 & 2a+2 & 2a+2 & 2-a & -3 
\end{bmatrix}\\
&= \underbrace{-\left[\begin{smallmatrix}
5 & 0 & 0 & 0 & 0\\
0 & 4 & 0 & 0 & 0\\
0 & 0 & 3 & 0 & 0\\
0 & 0 & 0 & 4 & 0\\
0 & 0 & 0 & 0 & 5
\end{smallmatrix}\right] + \left[\begin{smallmatrix*}[r]-3\\1\\-1\\-3\\2\end{smallmatrix*}\right]\left[\begin{smallmatrix}1\\1\\1\\1\\1\end{smallmatrix}\right]^{\T}}_{\define A_1} + a \left[\begin{smallmatrix}
2 & -1 & 0 & -1 & -1\\
-1 & 0 & -1 & 1 & -1\\
0 & 1 & -1 & 0 & 0\\
0 & 1 & -1 & 1 & 0\\
1 & 2 & 2 & -1 & 0
\end{smallmatrix}\right].
\end{aligned}
\]
We compute using \eqref{Eq:alpha} that $\mu = 0.35$, and hence $A_1$ is diagonally stable with certificate $D = \ddiag(\frac{1}{3}, 1,1,\frac{1}{3},\frac{1}{2})$. By applying \eqref{Eq:SigmaBound}, we conclude that $A$ is diagonally stable {\tb (again with certificiate $D$)} for all $|a|<0.1278$. \hfill \oprocend
\end{example}

\bigskip{}

\begin{remark}
Suppose that the (higher rank) interconnection matrix $S$ admits the singular value decomposition (SVD) 
\begin{equation}\label{Eq:GammaSVD}
S= \sum_{j=1}^N \nolimits \sigma_j u_j v_j^\T
\end{equation}
with a strictly dominant largest singular value, i.e, $$\sigma_1>\sigma_2\geq \ldots \geq \sigma_N.$$
We can decompose $S$ as
\[
S=\sigma_1 u_1v_1^T + E, \quad E:=S-\sigma_1u_1v_1^\top. 
\]
The first term on the right hand side of the above equality is the optimal rank-1 approximation of $S$ and the remainder matrix $E$ satisfies $\|E\|_2=\sigma_2$; this  is  the  so-called  Mirsky-Eckart-Young Theorem. Now if $v_1$ is a positive vector and $u_1$ does not contain any zero elements,\footnote{If instead $u_1\in \real_{>0}^N, \; v_{1} \in \real_{\neq 0}^{N}$ then one can obtain a similar condition by using the SVD of $S^\T$.}
then the treatment preceding \eqref{Eq:SigmaBound} yields the following condition for diagonal stability of $A=-\Delta+S$:
\begin{equation}\label{Eq:Condition-revised}
{\sigma_2} < \rho \left( 1-\sigma_1\sum_{i\in \mathcal{I}_N} \nolimits \frac{1}{\delta_i}[u_{1,i}v_{1,i}]_{+} \right),
\end{equation}
where 
\begin{align*}
&\rho(u_1,v_1,\Delta)=
\ddfrac{\min_{i\in \mathcal{I}_N}\,\delta_i \,(v_{1,i}/|u_{1,i}|)} {\max_{i\in \mathcal{I}_N}\,(v_{1,i}/|u_{1,i}|)}.
\end{align*}
The condition \eqref{Eq:Condition-revised} for diagonal stability states that (i) system $-\Delta$ coupled with \emph{only} the dominant mode $\sigma_1 u_1v_1^{\T}$ should be diagonally stable, and (ii) the norm of the approximation error matrix $E$, i.e, $\sigma_2$, should be sufficiently small. We note that the perturbed condition number $\rho$ is independent of the error matrix $E$, and depends only on $\Delta$ and on the dominant interconnection mode. A notable class of interconnection matrices with strictly dominant singular value $\sigma_1>\sigma_2$ and $v_1\in\real^N_{>0}$ (respectively, $u_1\in\real^N_{>0}$) is given by $S^\T S$ (respectively, $SS^\T$) being an irreducible nonnegative matrix.\footnote{A matrix $M \in \real^{N \times N}$ is called \textit{nonnegative} if all its elements are nonnegative, and is called \textit{irreducible} if the associated directed graph induced by $M$ is strongly connected; see \cite{RAH-CRJ:12} for details.} \hfill \oprocend
\end{remark}


\section{Stability of Automatic Generation Control in Power Systems}
\label{Sec:AGC}

We now leverage Theorem \ref{Thm:Main} as the key step in a larger stability treatment of automatic generation control in interconnected power systems. We begin in Section \ref{Sec:BackgroundAGC} by providing introductory background information on AGC. Section \ref{Sec:PowerModel} lays out our technical assumptions on the power system model under consideration, with Section \ref{Sec:ACE} formally introducing the AGC controller. In Section \ref{Sec:AGCStability} we state and prove the main stability result.

\subsection{Background on Automatic Generation Control}
\label{Sec:BackgroundAGC}

Large-scale AC power systems consist of interconnections between autonomously managed areas. Mismatch between generation and load within each area is compensated through a hierarchy of control layers operating on different spatiotemporal scales. The lowest ``primary'' control layer is decentralized, and uses proportional feedback of local frequency deviation to stabilize the system by quickly adjusting power generation levels; this occurs on a time-scale of seconds. The highest ``tertiary'' control layer is concerned with computing economically efficient generation set-points via global constrained optimization, through either economic dispatch \cite{AJW-BFW:96} or through the more sophisticated optimal power flow \cite{AJW-BFW:96,SHL:14a}; this occurs on a time-scale of 5 to 10 minutes. {\tb Our focus here is on the ``secondary'' layer of control, called Automatic Generation Control (AGC), which acts as a dynamic interface between the primary and tertiary control layers. The AGC system receives the optimal set-points from the tertiary control layer, and transmits modified set-points to the local primary controllers. The principal effect is that AGC eliminates generation-load mismatch \emph{within} each area, which in turn ensures that the system frequency and all net inter-area power exchanges are regulated to their scheduled values.}





Successfully deployed since the late 1940's \cite{NC:83}, AGC has a long and extensive history of study, and its evolution continues to be a topic of academic and industrial interest. We make no attempt to summarize the historical literature here. Standard textbook treatments can be found in \cite{AJW-BFW:96,PK:94}, and  \cite{FPD-RJM-WFB:73a,FPD-RJM-WFB:73b,HGK-KCK-AB:75,WBG:78,HG-JS:80,JC:85,MSC:86,NJ-LSV-DNE-LHF-AGH:92,IEEE-Report-AGC:79} contain outstanding historical and/or practitioner perspectives on AGC. An adjacent line of research considers the application of modern control techniques to secondary control. Several surveys \cite{PKI-DPK:05,HHA-MEHG-RZ-EHF-PS:18,DKM-FD-HS-SHL-SC-RB-JL:17,FD-SB-JWSP-SG:19a} are available summarizing aspects of this line of research, including distributed frequency control mechanisms.\footnote{It was suggested as early as 1978 \cite{ED-NP:78} that advanced decentralized control was unlikely to offer significant performance advantages {\tb in practice over traditional AGC. This conclusion has stood the test of time.}} {\tb Presently, heavy renewable energy integration is placing increased regulation demands on traditional AGC systems. As a result, several lines of work have arisen which extend traditional AGC schemes by modelling additional system dynamics in the design phase, or by explicitly accounting for sources of uncertainty.} Work in this direction includes model-predictive AGC \cite{ANV-IAH-JBR-SJW:08,AM-LI-KU:16,PRBM-PT:17}, various ``enhanced'' versions of AGC \cite{QL-MDI:12,DA:14,CL-JHC:17}, online gradient-type methods \cite{NL-CZ-LC:16, FD-SG:17, CZ-EM-SHL-JB:18, JWSP:20b}, and frequency or dynamics-aware dispatch and AGC \cite{DA-PWS-ADD:16,AAT-FZ-LX:11,GZ-JM-QW:19,EE-ZT-JWSP-EF-MP-HH:20k,PC-SD-YCC-MP:20}. While our subsequent analysis will be focused on the traditional AGC controller, our technical approach will likely prove useful for studying many of these variants, and for analyzing distributed consensus-based frequency control mechanisms (e.g., \cite{DKM-FD-HS-SHL-SC-RB-JL:17,FD-SB-JWSP-SG:19a}) in multi-area systems. 

Industry implementations of AGC are incredibly varied, and often include logical subroutines for reducing the activity of the AGC system and for handling system-specific conditions and operator preferences (e.g., \cite{FPD-RJM-WFB:73b,CWT-RLC:76}). The fundamental underlying characteristic that enables the success of these diverse implementations is the slow {\tb time constant of AGC} (minutes) relative to primary control dynamics (seconds). This slow speed is \emph{intentional}, as the goal of AGC is smooth re-balancing of generation and load inbetween economic re-dispatch. The slow speed is also \emph{necessary}: dynamic models of primary frequency dynamics (including energy conversion, turbine-governor, and load dynamics) are subject to considerable uncertainty, and sampling/communication/filtering processes introduce unavoidable delays. To ensure closed-loop stability, all practical AGC systems must be sufficiently slow so that no significant dynamic interaction occurs between the secondary loop and the primary frequency dynamics. As put succinctly in \cite{NJ-LSV-DNE-LHF-AGH:92}, ``\emph{reduction in the response time of AGC is neither possible nor desired}'' and ``\emph{attempting to do so serves no particular economic or control purposes.}'' This time-scale separation is exploited heavily in our subsequent analysis.

\subsection{Interconnected Power System Model}
\label{Sec:PowerModel}
We consider an interconnected power system consisting of $N$ areas, and label the set of areas as $\mathcal{A} = \{1,\ldots,N\}$. In area $k \in \mathcal{A}$, suppose there are $m_k$ generators; we label the set of generators as $\mathcal{G}_{k} = \{1,\ldots,m_k\}$, and we let $\mathcal{G}_{k}^{\rm AGC} \subseteq \mathcal{G}_k$ denote the subset of generators which participate in AGC. Each generator $i \in \mathcal{G}_{k}$ has a power-regulating turbine-governor system which accepts a power reference command $u_{k,i}$, with $u_{k,i}^{\star}$ denoting the base reference determined by the tertiary control layer. Power reference commands are restricted to the limits $\underline{u}_{k,i} \leq u_{k,i} \leq \overline{u}_{k,i}$, and we let
\[
\mathcal{U}_k \define \prod_{i\in\mathcal{G}_k} (\underline{u}_{k,i},\overline{u}_{k,i}), \qquad \mathcal{U} \define \prod_{k\in\mathcal{A}}\mathcal{U}_k
\]
denote the reference constraint sets. For generators not participating in AGC, we always have $u_{k,i} = u_{k,i}^{\star} = \underline{u}_{k,i} = \overline{u}_{k,i}$, and for notational convenience, we let
\begin{equation}\label{Eq:Deltau}
\Delta u_{k} \define  \sum_{j\in\mathcal{G}_{k}} (u_{k,j}-u_{k,j}^{\star}) = \sum_{j\in\mathcal{G}_{k}^{\rm AGC}} (u_{k,j}-u_{k,j}^{\star})
\end{equation}
denote the total change made in set-points for generation in area $k$ relative to the economic dispatch point. For area $k$ we define the vector variable $u_k = \mathrm{col}(u_{k,1},\ldots,u_{k,m_k})$, with $u_k^{\star}$ defined similarly. 

The \emph{net interchange} $\NI_k \in \real$ for area $k \in \mathcal{A}$ is defined as the net power transfer \emph{from} area $k$ \emph{to} the remainder of the interconnected system. Following the North American Electric Reliability Corporation (NERC) regulations \cite{NERC:10bal}, the power flow on any interconnecting tie line is measured by both areas at a common point, which implies that $\NI_k$ is simply the algebraic sum of power flows on all interconnecting tie lines. {\tb We let $f_k \in \real$ be a measurement of the AC frequency in area $k$, with $\Delta f_k = f_k - f_k^{\star}$ denoting the deviation from the scheduled value $f_k^{\star}$, and similarly define the net interchange deviation $\Delta \NI_{k} = \NI_{k} - \NI_{k}^{\star}$.} We collect all variables for \emph{all} areas into stacked vectors $u = \mathrm{col}(u_1,\ldots,u_N)$, with $\Delta f \in \real^N$ and $\Delta \NI \in \real^N$ similarly defined. The entire interconnected power system (all areas together, including primary stabilizing {\tb droop} controllers) is assumed to be described by the ODE model
{\tb
\begin{subequations}\label{Eq:NonlinearPowerSystem}
\begin{align}
\label{Eq:NonlinearPowerSystem-1}
\dot{x}(t) &= F(x(t),u(t),w)\\
\label{Eq:NonlinearPowerSystem-2}
\mathrm{col}(\Delta f(t),\Delta \NI(t)) &= h(x(t),u(t),w),
\end{align}
\end{subequations}
}
where $x(t) \in \tb{\mathcal{X}\subseteq \,} \real^n$ is the vector of states and $F,h$ are appropriate functions. {\tb The dynamics \eqref{Eq:NonlinearPowerSystem-1} are assumed to already include the typical low-pass filters for the measurements specified in \eqref{Eq:NonlinearPowerSystem-2}. The exogenous disturbance $w \in \real^{n_w}$ is assumed to be  constant, and can model unmeasured disturbances, set-points to the frequency/net interchange schedules, and references to other control loops.} Most importantly, $w$ includes the unmeasured net load deviation $\Delta P_{k}^{\rm L}$ for each area $k \in \mathcal{A}$ (see \eqref{Eq:SteadyStateFormulas}).

The precise model \eqref{Eq:NonlinearPowerSystem} is never known in practice, nor does the design of AGC systems rely on such an explicit model. Based on this engineering practice, we make no attempt to specify a detailed model, but will instead assume that the model satisfies some basic {\tb regularity,} stability, and steady-state properties. 


\smallskip

\begin{assumption}[\bf Interconnected Power System Model]\label{Ass:PowerSystem}
{There exist domains $\mathcal{X} \subseteq \real^n$ and $\mathcal{W} \subseteq \real^{n_w}$ such that the following hold:}
%
\begin{enumerate}
\item \label{Ass:PowerSystem-1} \textbf{Model Regularity:} $F$, $h$, and all associated Jacobian matrices are Lipschitz continuous on $\mathcal{X}$ uniformly in $(u,w) \in \mathcal{U} \times \mathcal{W}$;
\item \label{Ass:PowerSystem-2} \textbf{Existence of Steady-State:} there exists a continuously differentiable equilibrium map $\map{x_{\rm ss}}{\mathcal{U}\times\mathcal{W}}{\mathcal{X}}$ which is Lipschitz continuous on $\mathcal{U}\times\mathcal{W}$ and satisfies $0 = F(x_{\rm ss}(u,w),u,w)$ for all $(u,w) \in \mathcal{U} \times \mathcal{W}$;
\item \label{Ass:PowerSystem-3} \textbf{Uniform Exponential Stability of the Steady-State:} the steady-state $x_{\rm ss}(u,w)$ is locally exponentially stable, uniformly in the inputs $(u,w) \in \mathcal{U} \times \mathcal{W}$;
%
%
\item \label{Ass:PowerSystem-4} \textbf{Steady-State Synchronism and Interchange Balance:} for each $(u,w) \in \mathcal{U}\times\mathcal{W}$ the steady-state values of $(\Delta f,\Delta \NI)$ determined by $\mathrm{col}(\Delta f,\Delta \NI) = h(x_{\rm ss}(u,w),u,w)$ satisfy the synchronism condition
\begin{equation}\label{Eq:Synchronous}
\Delta f_1 = \Delta f_2 = \cdots = \Delta f_N
\end{equation}
and the net interchange balance condition
\begin{equation}\label{Eq:NetBalance}
0 = \sum_{k\in\mathcal{A}}\nolimits \Delta \NI_k.
\end{equation}
\item \label{Ass:PowerSystem-5} \textbf{Area Power Balance:}
the steady-state values $(\Delta f,\Delta \NI)$ from \ref{Ass:PowerSystem-4} satisfy the area-wise balance conditions
\begin{subequations}\label{Eq:SteadyStateFormulas}
\begin{align}
\label{Eq:Balancek}
\sum_{i\in\mathcal{G}_k}\nolimits (P_{k,i}-u_{k,i}^{\star}) &= \Delta P_{k}^{\rm FDL} + \Delta P_{k}^{\rm L} + \Delta\NI_{k}\\
\label{Eq:GovPower}
P_{k,i} &=  u_{k,i} - \tfrac{1}{R_{k,i}}\Delta f_k\\
\label{Eq:FDLPower}
\Delta P_{k}^{\rm FDL} &= D_k\Delta f_k
\end{align}
\end{subequations}
for each $k \in \mathcal{A}$ and $i \in \mathcal{G}_{k}$, where $P_{k,i}$ and $R_{k,i} > 0$ are the steady-state electrical power output and primary controller gain of generator $i \in \mathcal{G}_{k}$, and $\Delta P_{k}^{\rm FDL}$ and $D_{k} > 0$ are the steady-state power change due to frequency-dependent loads and load damping parameter of area $k$.
\end{enumerate}
\end{assumption}

\medskip

{\tb Assumptions \ref{Ass:PowerSystem-1}--\ref{Ass:PowerSystem-3} say that the model is sufficiently regular, and that for constant inputs $(u,w) \in \mathcal{U}\times\mathcal{W}$, the system possesses an exponentially stable equilibrium state $x = x_{\rm ss}(u,w) \in \mathcal{X}$. Note that this equilibrium need only be unique \emph{within} the set $\mathcal{X}$, which itself should be thought of as being contained in the normal operating region of the state space. While other equilibrium points may exist, they will typically be outside the normal operating range and therefore outside the set $\mathcal{X}$; see, e.g., \cite{JWSP:17b} for further discussion of non-uniqueness of equilibria in power system models.} Regarding \eqref{Eq:Synchronous}, all AC power systems self-synchronize under normal operating conditions (see, e.g., \cite{fd-fb:09z}), and we restrict our attention to such cases. The net interchange balance condition \eqref{Eq:NetBalance} holds by definition of the net interchange, as any power leaving one area must necessarily enter another. Similarly, the equality \eqref{Eq:Balancek} describes the balance of power for each area, with change in injected power on the left balancing change in extracted power on the right, which consists of net load change (plus incremental losses) $\Delta P_{k}^{\rm L}$,  net interchange deviation $\Delta \NI_k$, and frequency-dependent load effects $\Delta P_{k}^{\rm FDL}$. The preceding assumptions all arise from the fundamental physics and engineering design of AC power systems. Our core simplifying modelling assumptions are \eqref{Eq:GovPower} and \eqref{Eq:FDLPower}, which state that {\tb generator droop controls and any frequency-dependent loads} respond \emph{linearly} to deviations in frequency. These assumptions are standard in textbook equilibrium analysis of AGC \cite{PK:94}.

{\tb To avoid a further analysis of saturated operation and integrator anti-windup implementations for AGC, we assume there is sufficient regulation capacity in each area.}

\smallskip

\begin{assumption}[\bf Strict Local Feasibility]\label{Ass:Feasibility}
Each area has sufficient regulation capacity to meet the disturbance, i.e.,
\[
\Delta P^{\rm L}_{k} \in \mathcal{C}_k \define \Bigg(\sum_{i\in\mathcal{G}_k^{\rm AGC}} (\underline{u}_{k,i}-u_{k,i}^{\star}), \sum_{i\in\mathcal{G}_k^{\rm AGC}}(\overline{u}_{k,i}-u_{k,i}^{\star})\Bigg)
\]
for each area $k \in \mathcal{A}$.
\end{assumption}

\smallskip

{\tb
\begin{remark}[\bf Power System Models]\label{Rem:ODE}
Power system models are usually expressed as differential-algebraic equations (DAEs), and the ODE model \eqref{Eq:NonlinearPowerSystem} may have been obtained from a DAE through a number of techniques. Engineering techniques for obtaining ODE models from DAE power system models include Kron reduction of load buses \cite{FD-FB:11d} and the introduction of frequency-dependent loads \cite{ARB-DJH:81}. More broadly, DAE power system models always occur in the so-called semi-explicit form
\begin{equation}\label{Eq:DAE}
\dot{x}_1 = g_1(x_1,x_2,u,w), \qquad 0 = g_2(x_1,x_2,u,w)
\end{equation}
where $x_1$ are the differential variables and $x_2$ are the algebraic variables. When $\frac{\partial g_2}{\partial x_2}$ is non-singular for all $(x_1,x_2,u,w)$ on a domain of interest, the DAE \eqref{Eq:DAE} is said to have index-1, and can be reduced to an ODE model of the form \eqref{Eq:NonlinearPowerSystem}. This reduction can sometimes be done explicitly via elimination and/or differentiation (e.g., \cite{NM-CDP-AJvdS-JMAS:18}), but more generally is an implicit reduction via the inverse function theorem \cite{DJH-IMYM-90}. As we do not require an explicit expression for the functions $F$ and $h$ in \eqref{Eq:NonlinearPowerSystem}, such an implicit reduction poses no problems for our analysis.
\hfill \oprocend
\end{remark}
}

\smallskip

%
%
%

\subsection{Area Control Error and AGC Model}
\label{Sec:ACE}

As discussed in Section \ref{Sec:BackgroundAGC}, the purpose of AGC is to asymptotically eliminate generation-load mismatch within each balancing area. Unfortunately, the net load change $\Delta P_{k}^{\rm L}$ is exogenous and typically unmeasurable; hence, even if all generation was metered, the generation-load mismatch cannot be computed for use in a regulating controller. Standard AGC implementations circumvent this issue by defining an auxiliary error signal called the \emph{area control error} $\mathsf{ACE}_k$, which is defined by taking the net-interchange deviation $\Delta \NI_k$ and \emph{biasing} it using the measured frequency deviation
\begin{equation}\label{Eq:ACE}
\begin{aligned}
\mathsf{ACE}_{k}(t) &\define \Delta \NI_{k}(t) + b_k\Delta f_k(t),\\
\end{aligned}
\end{equation}
where $b_k > 0$ is the \emph{frequency bias} for area $k$. Surprisingly, each area individually zeroing this error signal achieves the desired re-balancing of generation and demand.

\smallskip

\begin{lemma}[\bf Steady-State Zeroing of ACE]\label{Lem:ACE}
If the interconnected power system \eqref{Eq:NonlinearPowerSystem} is in steady-state as specified in Assumption \ref{Ass:PowerSystem}, then the following statements are equivalent: 
\begin{enumerate}
\item \label{Lem:ACE-1} $\Delta u_k = \Delta P_{k}^{\rm L}$ for all areas $k \in \mathcal{A}$;
\item \label{Lem:ACE-2} $\Delta f_{k} = 0$ and $\Delta \NI_{k} = 0$ for all areas $k \in \mathcal{A}$;
\item \label{Lem:ACE-3} $\mathsf{ACE}_k = 0$ for all areas $k \in \mathcal{A}$.
\end{enumerate}
\end{lemma}

\smallskip

Put simply, balancing power in each area is equivalent to simultaneous regulation of frequency and all net-interchange deviations, which is in turn equivalent to zeroing of all area control errors.

\begin{pfof}{Lemma \ref{Lem:ACE}} Begin by noting that substitution of \eqref{Eq:GovPower},\eqref{Eq:FDLPower} into \eqref{Eq:Balancek} yields the compact equalities
\begin{equation}\label{Eq:ReducedBalanceEqn}
\Delta u_k = \Delta \NI_k + \beta_k\Delta f_k + \Delta P_{k}^{\rm L}, \quad k \in \mathcal{A},
\end{equation}
where $\beta_k \define D_k + \sum_{i\in\mathcal{G}_k} R_{k,i}^{-1}$.

\ref{Lem:ACE-1} $\Longleftrightarrow$ \ref{Lem:ACE-2}: For the forward direction, \eqref{Eq:ReducedBalanceEqn} implies that
\begin{equation}\label{Eq:ACEProof1}
\Delta\NI_k + \beta_k \Delta f_k = 0, \quad k \in \mathcal{A}.
\end{equation}
Summing \eqref{Eq:ACEProof1} over all areas $k$ and using \eqref{Eq:NetBalance}, we find that $0 = \sum_{k\in\mathcal{A}}\beta_k\Delta f_k$ which using \eqref{Eq:Synchronous} further implies that $\Delta f_{k} = 0$ for all $k \in \mathcal{A}$. It follows now from \eqref{Eq:ACEProof1} that $\Delta \NI_k = 0$ for all $k$ as well. The converse is immediate by setting $\Delta f_k = \Delta \NI_k = 0$ in \eqref{Eq:ReducedBalanceEqn}.

\ref{Lem:ACE-2} $\Longleftrightarrow$ \ref{Lem:ACE-3}: The forward direction is immediate from \eqref{Eq:ACE}. The argument for the converse implication is identical to that following \eqref{Eq:ACEProof1} and is omitted.
\end{pfof}

\begin{figure}[ht!]
\centering
\includegraphics[width=1\columnwidth]{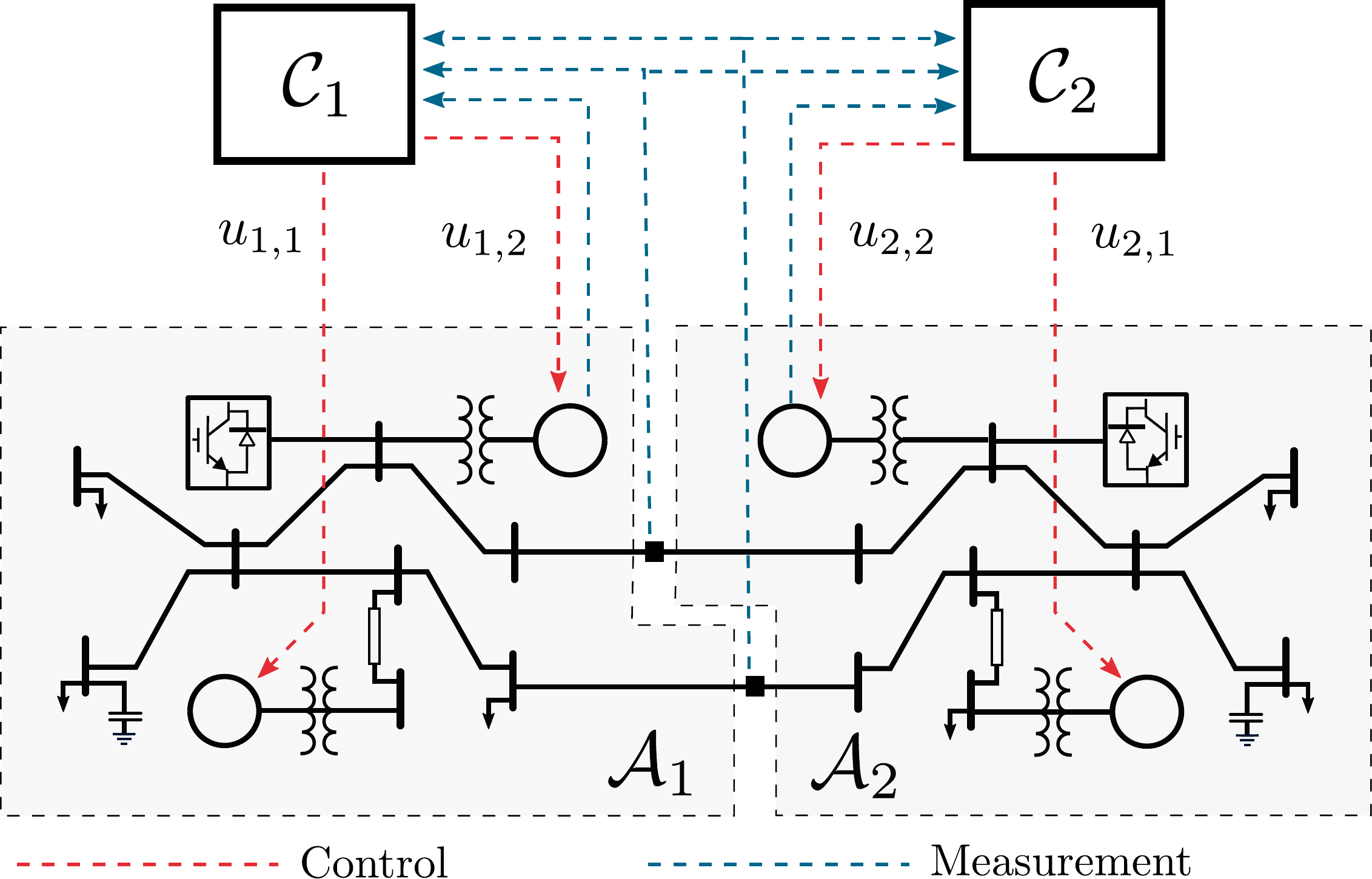}
\caption{Diagram of AGC in a two-area interconnected power system. The power flows on both interconnecting tie lines are measured at common points, while frequency measurements are taken locally from a chosen generator within each area.}
\label{Fig:AGC}
\end{figure}

Each area $k \in \mathcal{A}$ integrates the ACE to produce a local AGC control signal $\eta_k$ as
\begin{equation}\label{Eq:AGCSimple}
\tau_k \dot{\eta}_k(t) = -\mathsf{ACE}_k(t), \qquad k \in \mathcal{A},
\end{equation}
where $\tau_k > 0$ is the integral time constant, quoted in the literature as ranging from 30s up to 200s. 
Control actions from the AGC system are allocated across all participating generators $\mathcal{G}_{k}^{\rm AGC}$ such that their incremental costs of production (or, with lossess, delivery) are roughly equalized \cite{AJW-BFW:96}. A typical allocation rule including limiting of control signals is
\begin{equation}\label{Eq:Allocation}
u_{k,i} = \mathrm{sat}_{k,i}(u_{k,i}^{\star} + \alpha_{k,i}\eta_k), \quad i \in \mathcal{G}_k^{\rm AGC},
\end{equation}
where 
\[
\mathrm{sat}_{k,i}(v) = \begin{cases}
v &\text{if}\,\, v \in (\underline{u}_{k,i},\overline{u}_{k,i})\\
\overline{u}_{k,i} &\text{if}\,\, v > \overline{u}_{k,i}\\
\underline{u}_{k,i} &\text{if}\,\, v < \underline{u}_{k,i}
\end{cases}
\]
and $\{\alpha_{k,i}\}_{i\in\mathcal{G}_{k}^{\rm AGC}}$ are nonnegative \emph{participation factors} \cite{AJW-BFW:96} with normalization
\[
\sum_{i\in\mathcal{G}_{k}^{\rm AGC}}\nolimits \alpha_{k,i} = 1, \quad k \in \mathcal{A}.
\]
A schematic of AGC in a two-area power system is shown in Figure \ref{Fig:AGC}; note that the control scheme is area-wise decentralized, and no communication occurs between the two areas.

\smallskip


\subsection{Closed-Loop Stability with AGC}
\label{Sec:AGCStability}

The closed-loop system consists of the interconnected power system model \eqref{Eq:NonlinearPowerSystem} with the decentralized AGC controllers \eqref{Eq:AGCSimple} and the allocation rules \eqref{Eq:Allocation}. We now proceed with a stability analysis based on time-scale separation, where we assume the AGC dynamics are slow compared to the power system dynamics \eqref{Eq:NonlinearPowerSystem}. This analysis approach is strongly justified by the time-scale properties of AGC (Section \ref{Sec:BackgroundAGC}); the slowest time-constant in \eqref{Eq:NonlinearPowerSystem} will be that of the primary control dynamics, which is on the order of seconds, while the time-constant of AGC is on the order of minutes. We can now state the main result.


\smallskip

\begin{theorem}[\bf Stability with AGC]\label{Thm:AGCStable}    
Consider the interconnected power system \eqref{Eq:NonlinearPowerSystem} under Assumptions \ref{Ass:PowerSystem}--\ref{Ass:Feasibility} with the AGC controllers \eqref{Eq:AGCSimple} and the allocation rules \eqref{Eq:Allocation}. There exists a value $\tau^{\star} > 0$ such that if $\min_{k\in\mathcal{A}}\tau_k \geq \tau^{\star}$, then
\begin{enumerate}
\item[(a)] the closed-loop system possesses a unique exponentially stable equilibrium point $(\bar{x},\bar{\eta}) \in \mathcal{X} \times \real^{N}$, and
\item[(b)] $\mathsf{ACE}_k(t) \to 0$ as $t \to \infty$ for all areas $k \in \mathcal{A}$.
\end{enumerate}
\end{theorem}

\smallskip

Theorem \ref{Thm:AGCStable} states that \textemdash{} with the usual time-scales of operation \textemdash{} closed-loop stability with AGC is guaranteed for \emph{any} tuning of bias factors $b_k > 0$ and with no assumptions on homogeneity of time constants $\tau_k$. This rather ``unconditional'' stability result is consistent with engineering practice, in which balancing authorities independently tune their AGC controllers without coordination. While an explicit expression for $\tau^{\star}$ can in fact be obtained under our assumptions, the information required to compute it would not be available in practice, and hence the result is most meaningfully stated in the qualitative form above.

\smallskip

\begin{pfof}{Theorem \ref{Thm:AGCStable}} {\tb Let $\tau > 0$, and for each $k \in \mathcal{A}$ define $\tilde{\tau}_k \define \tau_k/\tau$. With this, the time constants in \eqref{Eq:AGCSimple} can be expressed as $\tau_k = \tau \cdot \tilde{\tau}_k$. Set $\varepsilon \define 1/\tau$ and define the new time variable $\ell \define \varepsilon t$. With this, the closed-loop system \eqref{Eq:NonlinearPowerSystem},\eqref{Eq:AGCSimple},\eqref{Eq:Allocation} can be written in standard singularly perturbed form \cite{AS-HKK:84},\cite[Chp. 11]{HKK:02} as}
\begin{subequations}\label{Eq:SingularPerturbation}
\begin{align}
\label{Eq:SingularPerturbation-1}
\varepsilon \frac{\mathrm{d}x}{\mathrm{d}\ell} &= F(x,u,w)\\
\label{Eq:SingularPerturbation-1b}
\mathrm{col}(\Delta f,\Delta \NI) &= h(x,u,w)\\
\label{Eq:SingularPerturbation-2}
\tilde{\tau}_k\frac{\mathrm{d}\eta_k}{\mathrm{d}\ell} &= -(\Delta \NI_k + b_k\Delta f_k)\\
\label{Eq:SingularPerturbation-3}
u_{k,i} &= \mathrm{sat}_{k,i}(u_{k,i}^{\star} + \alpha_{k,i}\eta_k),
\end{align}
\end{subequations}
for $k \in \mathcal{A}$ and $i \in \mathcal{G}_{k}^{\rm AGC}$. Our approach will be to invoke the singular perturbation stability result \cite[Lemma 11.4]{HKK:02}, and we must therefore check the stipulated conditions on the boundary layer and reduced dynamics. The boundary layer dynamics are \eqref{Eq:SingularPerturbation-1} with $(u,w)$ considered as parameters, and Assumption \ref{Ass:PowerSystem} \ref{Ass:PowerSystem-1}--\ref{Ass:PowerSystem-3} are sufficient to meet the conditions on the boundary layer system imposed in \cite[Lemma 11.4]{HKK:02}.
%
%
To obtain the reduced dynamics associated with \eqref{Eq:SingularPerturbation}, we set $\varepsilon = 0$ in \eqref{Eq:SingularPerturbation-1} and use Assumption \ref{Ass:PowerSystem} \ref{Ass:PowerSystem-4} and \ref{Ass:PowerSystem-5}. First, substituting \eqref{Eq:GovPower} into \eqref{Eq:Balancek}, summing over $k \in \mathcal{A}$, and using \eqref{Eq:Synchronous} and \eqref{Eq:NetBalance}, one finds that
\begin{equation}
\label{Eq:SteadyStateFrequency}
\Delta f_{k} = \tfrac{1}{\beta}\sum_{i\in\mathcal{A}}\nolimits (\Delta u_i - \Delta P_{i}^{\rm L}), \quad k \in \mathcal{A},
\end{equation}
where $\Delta u_k$ is as defined in \eqref{Eq:Deltau}, with $\beta \define \sum_{k \in \mathcal{A}}\beta_k$ and $\beta_k \define D_k + \sum_{i\in\mathcal{G}_{k}}R_{k,i}^{-1}$. Substituting \eqref{Eq:SteadyStateFrequency} back into \eqref{Eq:Balancek}, one similarly obtains the steady-state relationship
\begin{equation}\label{Eq:SteadyStatePowerFlow}
\begin{aligned} 
\Delta \NI_{k} &= \tfrac{\beta-\beta_k}{\beta}(\Delta u_k - \Delta P_{k}^{\rm L})\\
&\qquad - \tfrac{\beta_k}{\beta}\sum_{j\in\mathcal{A}\setminus\{k\}}\nolimits (\Delta u_j - \Delta P_{j}^{\rm L}).
\end{aligned}
\end{equation}
Substituting \eqref{Eq:SteadyStateFrequency} and \eqref{Eq:SteadyStatePowerFlow} into \eqref{Eq:ACE}, we therefore find that in steady-state
\begin{equation}\label{Eq:SteadyStateACE}
\begin{aligned}
\mathsf{ACE}_{k} &= \tfrac{\beta+b_k-\beta_k}{\beta}(\Delta u_k - \Delta P_{k}^{\rm L})\\
&\qquad  + \tfrac{b_k-\beta_k}{\beta}\sum_{j \in \mathcal{A}\setminus\{k\}}\nolimits (\Delta u_j - \Delta P_{j}^{\rm L}).
\end{aligned}
\end{equation}
Using \eqref{Eq:SingularPerturbation-3}, note that we may compactly write
\begin{equation}\label{Eq:DeltaukSteady}
\begin{aligned}
\Delta u_k &=  \sum_{i\in\mathcal{G}_k}\nolimits(\mathrm{sat}_{k,i}(u_{k,i}^{\star} + \alpha_{k,i}\eta_k) - u_{k,i}^{\star})\\
&\define \varphi_{k}(\eta_k).
\end{aligned}
\end{equation}
Substituting \eqref{Eq:SteadyStateACE} and \eqref{Eq:DeltaukSteady} back into \eqref{Eq:SingularPerturbation-2} and writing everything in vector notation, we obtain the reduced AGC dynamics
\begin{equation}\label{Eq:SlowTimeScaleAGC2}
\tilde{\tau} \dot{\eta} = \mathcal{B}(\varphi(\eta)-\Delta P^{\rm L})
\end{equation}
where $\tilde{\tau} = \mathrm{diag}(\tilde{\tau}_1,\ldots,\tilde{\tau}_N)$, $\eta = \mathrm{col}(\eta_1,\ldots,\eta_N)$, $\Delta P^{\rm L} = \mathrm{col}(\Delta P^{\rm L}_{1},\ldots,\Delta P^{\rm L}_{N})$, and $\varphi(\eta) = \mathrm{col}(\varphi_1(\eta_1),\ldots,\varphi_{N}(\eta_N))$. The matrix $\mathcal{B}$ is given by
\begin{equation*}
\mathcal{B} \define -\frac{1}{\beta}
{\footnotesize
\begin{bmatrix}
\beta + b_1 - \beta_1 & b_1 - \beta_1 & \cdots & b_1 - \beta_1\\
b_2 - \beta_2 & \beta + b_2 - \beta_2 & \cdots & \cdots \\
\vdots & \cdots & \ddots & b_{N-1} - \beta_{N-1}\\
b_N - \beta_N & \cdots & b_N  - \beta_N & \beta + b_N - \beta_N 
\end{bmatrix},
}
\end{equation*}
which can be written compactly in vector notation as
\begin{equation}\label{Eq:BAGC}
\mathcal{B} = -I_N + \tfrac{1}{\beta}(\boldsymbol{\beta}-\boldsymbol{b})\vones[N]^{\T}.
\end{equation}
where $\boldsymbol{\beta} = \mathrm{col}(\beta_1,\ldots,\beta_N)$ and $\boldsymbol{b} = \mathrm{col}(b_1,\ldots,b_N)$. In particular, $\mathcal{B}$ satisfies the assumptions of Theorem \ref{Thm:Main} with $\Delta = I_N$, $x = \tfrac{1}{\beta}(\boldsymbol{\beta}-\boldsymbol{b})$ and $y = \vones[N]$. Moreover, since $b_k, \beta_k > 0$, we have that
\[
\sum_{k=1}^{N}\nolimits \frac{1}{\delta_k}[x_ky_k]_{+} = \tfrac{1}{\beta}\sum_{k=1}^{N}\nolimits[\beta_k-b_k]_{+} < \frac{\sum_{k=1}^{N}\beta_k}{\beta} = 1.
\]
We therefore conclude from Theorem \ref{Thm:Main} that $\mathcal{B}$ is diagonally stable, and we let $D = \mathrm{diag}(d_1,\ldots,d_N) \succ 0$ denote a certificate such that $Q \define -\tfrac{1}{2}(\mathcal{B}^{\T}D + D\mathcal{B}) \succ 0$.

We now argue that the reduced system \eqref{Eq:SlowTimeScaleAGC2} possesses a unique equilibrium point. By strict feasibility (Assumption \ref{Ass:Feasibility}), we have that $\Delta P^{\rm L}_{k} \in \mathcal{C}_k$. Examining \eqref{Eq:DeltaukSteady}, it is straightforward to argue that $\map{\varphi_k}{\real}{\real}$ is piecewise linear, non-decreasing, and that $\mathrm{image}(\varphi_k) = \mathrm{closure}(\mathcal{C}_k)$. The preimage $\mathcal{P}_k = \setdef{\eta_k \in \real}{\varphi_{k}(\eta_k) \in \mathcal{C}_k}$ of $\mathcal{C}_k$ under $\varphi_k$ can be explicitly computed to be the interval
\[
\mathcal{P}_k \define \left(\min_{i\in\mathcal{G}_k}\tfrac{1}{\alpha_{k,i}}(\underline{u}_{k,i}-u_{k,i}^{\star}),\max_{i\in\mathcal{G}_k}\tfrac{1}{\alpha_{k,i}}(\overline{u}_{k,i}-u_{k,i}^{\star})\right),
\]
and a simple argument shows that $\varphi_k$ is a \emph{strictly} increasing function on $\mathcal{P}_k$. We conclude that the restriction $\map{\varphi_{k}|_{\mathcal{P}_k}}{\mathcal{P}_k}{\mathcal{C}_k}$ is a bijective function, and hence for each $k \in \mathcal{A}$, there exists a unique $\bar{\eta}_k \in \mathcal{P}_k$ such that $\varphi_{k}(\bar{\eta}_k) = \Delta P^{\rm L}_{k}$. Since $\mathcal{B}$ is diagonally stable, it is invertible, and we conclude that $\bar{\eta} = \mathrm{col}(\bar{\eta}_1,\ldots,\bar{\eta}_N)$ is the unique equilibrium point of the reduced dynamics \eqref{Eq:SlowTimeScaleAGC2}. 

For this equilibrium point we now define a Lyapunov candidate $\map{V}{\real^N}{\real}$ as
\begin{equation}\label{Eq:Popov}
V(\eta) = \sum_{k=1}^{N} d_k \tilde{\tau}_k\int_{\bar{\eta}_k}^{\eta_{k}} (\varphi_k(\xi_k) - \varphi_k(\bar{\eta}_k))\,\mathrm{d}\xi_k.
\end{equation}
Note that $V$ is continuously differentiable and that $V(\bar{\eta}) = 0$. Due to the previously noted properties of the function $\varphi_k$ on the open interval $\mathcal{P}_k$, there exist constants $\ell_k, L_k, r_k > 0$ such that
\[
\ell_k |\eta_k-\bar{\eta}_k|^2 \leq (\eta_k-\bar{\eta}_k)(\varphi_k(\eta_k) - \varphi_k(\bar{\eta}_k)) \leq L_{k} |\eta_k-\bar{\eta}_k|^2
\]
for all $\eta_k \in [\bar{\eta}_k-r_k,\bar{\eta}_k+r_k]$. It follows quickly from \eqref{Eq:Popov} then that there exist constants $c_5, c_6, c_7, r > 0$ such that
\begin{subequations}\label{Eq:LyapReduced}
\begin{align}
\label{Eq:LyapReduced-a}
c_5 \|\eta-\bar{\eta}\|_2^2 \leq V(\eta) &\leq c_6 \|\eta-\bar{\eta}\|_2^2\\
\label{Eq:LyapReduced-b}
\|\nabla V(\eta)\|_2 &\leq c_7 \|\eta-\bar{\eta}\|_2
\end{align}
\end{subequations}
for all $\eta$ such that $\|\eta-\bar{\eta}\|_2 < r$. We can now compute along trajectories of \eqref{Eq:SlowTimeScaleAGC2} that
\begin{equation}\label{Eq:LyapReducedDecrease}
\begin{aligned}
\dot{V}(\eta) &= (\varphi(\eta)-\varphi(\bar{\eta}))^{\T}D\mathcal{B}(\varphi(\eta)-\Delta P^{\rm L}_{k})\\
&= (\varphi(\eta)-\varphi(\bar{\eta}))^{\T}D\mathcal{B}(\varphi(\eta)-\varphi(\bar{\eta}))\\
&= \tfrac{1}{2}(\varphi(\eta)-\varphi(\bar{\eta}))^{\T}(D\mathcal{B}+\mathcal{B}^{\T}D)(\varphi(\eta)-\varphi(\bar{\eta}))\\
&= -(\varphi(\eta)-\varphi(\bar{\eta}))^{\T}Q(\varphi(\eta)-\varphi(\bar{\eta}))\\
&\leq -c_8 \|\eta-\bar{\eta}\|_2^2
\end{aligned}
\end{equation}
for some constant $c_8 > 0$ and all $\eta$ such that $\|\eta-\bar{\eta}\|_2 < r$. We conclude from \cite[Theorem 4.10]{HKK:02} that $\bar{\eta}$ is a locally exponentially stable equilibrium point for the reduced dynamics \eqref{Eq:SlowTimeScaleAGC2}, and moreover, that the conditions \eqref{Eq:LyapReduced} and \eqref{Eq:LyapReducedDecrease} meet the remaining requirements imposed in \cite[Lemma 11.4]{HKK:02}. {\tb It now follows that there exists $\varepsilon^{\star} > 0$ such that for all $\varepsilon \in (0,\varepsilon^{\star})$, the equilibrium $(x_{\rm ss}(\bar{u},w),\bar{\eta})$ of \eqref{Eq:SingularPerturbation} \textemdash{} and hence, of the closed-loop system \textemdash{} is locally exponentially stable. Statement (a) of the theorem is now immediately obtained by setting $\tau^{\star} = \frac{1}{\varepsilon^{\star}}\min_{k}\tilde{\tau}_k$. Finally, since $\mathsf{ACE} = - \mathcal{B}(\varphi(\bar{\eta})-\Delta P^{\rm L}) = 0$ at equilibrium, we conclude that $\mathsf{ACE}_k \to 0$ as $t \to \infty$, which shows statement (b) and completes the proof.} \end{pfof}

{\tb The result of Theorem \ref{Thm:AGCStable} can be interpreted as a rigorous dynamic systems justification of the experiential observation that low-gain AGC systems lead to stable interconnected power systems. We note that our main diagonal stability result Theorem \ref{Thm:Main} is essential in the proof of Theorem \ref{Thm:AGCStable}, as it allows us to use the composite Lyapunov construction \eqref{Eq:Popov} to show stability of the reduced dynamics \eqref{Eq:SlowTimeScaleAGC2}.
}

{\tb 
\subsection{Implications for Dynamic Performance of AGC}
\label{Sec:DynamicPerformance}

We now explore the implications of our analysis in Theorem \ref{Thm:AGCStable} for tuning and dynamic performance of AGC. Our starting point is the reduced dynamics \eqref{Eq:SlowTimeScaleAGC2} developed in the proof of Theorem \ref{Thm:AGCStable}. These dynamics describe the dynamics of AGC in the interconnected power system over a long time-scale (i.e., after the transient action of primary droop controllers). Further scrutiny of these reduced dynamics will reveal some of the fundamental performance characteristics of AGC, which arise due to the decentralized control structure and due to the selection of frequency bias gains $b_k$. 
%
%

\subsubsection{Ideal Bias Tuning for Non-Interaction}
\label{Sec:IdealTuning}

The ``optimal'' choice of the frequency bias gains has been a topic of substantial historical interest and controversy. In \cite{NC:56}, Cohn argued\footnote{The argument is based on a static equilibrium analysis; an accessible treatment can be found in \cite{PK:94}.} that each area should ideally set its bias $b_k$ equal to its \emph{frequency characteristic} $\beta_k = D_k + \sum_{i\in\mathcal{G}_k} R_{k,i}^{-1}$, and that in doing so, each area will minimally respond to disturbances occurring in other areas.

In the context of our reduced dynamics  \eqref{Eq:SlowTimeScaleAGC2}, if the bias $b_k$ is tuned such that $b_k = \beta_k$, then all off-diagonal elements in the $k$th row of $\mathcal{B}$ become zero. The $k$th equation in the dynamics \eqref{Eq:SlowTimeScaleAGC2} then decouples from the other states, and simplifies to the single-input single-output system
\begin{equation}\label{Eq:Decoupled}
\tilde{\tau}_k \dot{\eta}_k = -\varphi_{k}(\eta_k) + \Delta P^{\rm L}_{k}.
\end{equation}
This shows that the control variable $\eta_k$ for area $k$ converges with what is essentially a first-order response to the disturbance value $\Delta P_{k}^{\rm L}$, and is not influenced by the disturbances or control actions in any other areas. If \emph{all} areas select $b_k=\beta_k$, then all AGC systems are non-interacting; this provides a dynamic systems justification for Cohn's conclusion. While the tuning $b_k = \beta_k$ is recognized as the ideal one by regulatory bodies \textemdash{} including NERC in North America and ENTSO-E in Europe \cite{NERC:11,UCTE-App1:04}  \textemdash{} it unfortunately cannot be implemented in practice, as $\beta_k$ cannot be reliably estimated \cite{NERC:11}; this leads to our next study.

\smallskip

\subsubsection{Stability Margin with Overbiasing vs. Underbiasing}
\label{Sec:Overbiasing}

The tuning of the AGC system \eqref{Eq:AGCSimple} for area $k$ is said to be \emph{overbiased} if $b_k > \beta_k$, and is \emph{underbiased} if $b_k < \beta_k$. As $\beta_k$ cannot be reliably estimated in practice, NERC guidelines recommend that system operators err towards overbiasing the tunings of their AGC systems as opposed to underbiasing \cite{NERC:10bal}. In particular, in the US Eastern Interconnection, established procedures for setting $b_k$ based on peak load are thought to lead to a  100\% over-biasing. Our Lyapunov analysis of the reduced dynamics \eqref{Eq:SlowTimeScaleAGC2} can be used to generate insights into the effects of overbiasing vs. underbiasing.

To examine the aggregate effect of overbiasing and underbiasing, we consider a uniform tuning setup where each bias $b_k$ is directly proportional to the frequency characteristic $\beta_k$, expressed as $b_k = \kappa \beta_k$ for each $k \in \mathcal{A}$ and some $\kappa>0$. The case $\kappa \in (0,1)$ would be a (globally) underbiased tuning, while $\kappa \in (1,\infty)$ would be a (globally) overbiased tuning. Assuming $\kappa \neq 1$, we can leverage Corollary \ref{Cor:Explicit} to develop an explicit expression for the coefficients $d_k$ used in the Lyapunov construction \eqref{Eq:Popov}. In particular, after rescaling one concludes that $d_k = \tfrac{\beta}{\beta_k}$ is an eligible selection, where $\beta = \sum_{k\in\mathcal{A}}\beta_k$. With $D = \ddiag(d_1,\ldots,d_N)$ we can now compute that
\[
\begin{aligned}
D\mathcal{B} &= \beta\ddiag(1/\beta_k)\left(-I_N + \tfrac{1}{\beta}(\boldsymbol{\beta}-\boldsymbol{b})\vones[N]^{\T}\right)\\
&= -\beta\ddiag(1/\beta_k) + (1-\kappa)\vones[N]\vones[N]^{\T},
\end{aligned}
\]
and therefore
\[
\begin{aligned}
Q = -\tfrac{1}{2}(D\mathcal{B}+\mathcal{B}^{\T}D) &= \beta\ddiag(1/\beta_k) - (1-\kappa)\vones[N]\vones[N]^{\T}\\
&= \beta\ddiag(1/\beta_k)\left(I_N - \tfrac{1-\kappa}{\beta}\boldsymbol{\beta}\vones[N]^{\T}\right).
\end{aligned}
\]
From \eqref{Eq:LyapReducedDecrease}, the minimum eigenvalue of ${Q}$ provides a bound on the stability margin of the reduced AGC dynamics. To probe further, note that the matrix in brackets above has $N-1$ eigenvalues at $1$, with $N$th eigenvalue equal to $\kappa$, and hence
\[
\lambda_{\rm min}(Q) \geq \tfrac{\beta}{\min_{k\in\mathcal{A}}\beta_k}\cdot \min\{\kappa,1\}.
\]
We conclude that underbiased tunings ($\kappa < 1$) degrade the guaranteed margin of exponential stability, while overbiased tunings ($\kappa > 1$) do not. This provides a system-theoretic interpretation of NERC's preference for overbiased tunings vis-a-vis stability.


\smallskip

\subsubsection{Response of Area Control Errors to Disturbances}
\label{Sec:ResponseToDisturbances}

For our final study we will use the reduced dynamics \eqref{Eq:SlowTimeScaleAGC2} to examine the response of the area control errors $\mathsf{ACE}_k$ (i.e., the regulated variable) to changes in the load disturbances $\Delta P^{\rm L}_{k}$. To leverage LTI analysis tools, we will ignore the effects of saturation, so that $\varphi(\eta) = \eta$. Moreover, to focus in particular on the effects of bias tuning, we assume equal time constants  $\tilde{\tau}_k = \tilde{\tau}_{\ell} = \tau^{\prime}$ for some $\tau^{\prime} > 0$ and all $k,\ell \in \mathcal{A}$. Under these assumptions \eqref{Eq:SlowTimeScaleAGC2} \textemdash{} with ACE values as outputs of interest \textemdash{} becomes the LTI system
\begin{equation}\label{Eq:SlowTimeScaleAGC2-LTI}
\begin{aligned}
 \dot{\eta} &= -\tfrac{1}{\tau^{\prime}}\mathcal{B}(\eta-\Delta P^{\rm L}), \quad \mathsf{ACE} = \mathcal{B}(\eta-\Delta P^{\rm L}).
 \end{aligned}
\end{equation}
The transfer functions of interest can now be defined as
\[
S_{ij}(s) = \frac{\mathsf{ACE}_i(s)}{\Delta P^{\rm L}_{j}(s)} = -\mathsf{e}_i^{\T}(\tau^{\prime} sI_{N}-\mathcal{B})^{-1}\mathcal{B}\mathsf{e}_j,
\]
where $\mathrm{e}_i$ is the $i$th unit vector of $\real^N$ and where, in the second equality, we have substituted \eqref{Eq:SlowTimeScaleAGC2-LTI} and simplified. Through straightforward but tedious calculation, we find that
\[
\begin{aligned}
S_{ij}(s) 
&= -\frac{\tau^{\prime}s}{\tau^{\prime}s+1}\left[\delta_{ij} - \frac{1}{\beta}(\beta_i-b_i)\frac{\tau^{\prime}s}{\tau^{\prime}s+\tfrac{1}{\beta}\sum_{k\in\mathcal{A}}b_k}\right],
\end{aligned}
\]
where $\delta_{ij}$ is the usual Kronecker delta function. A representative Bode plot of $S_{ii}(s)$ for underbiased and overbiased tunings is shown in Figure \ref{Fig:Bode}. The peak sensitivity of this transfer function can easily be computed to be
\begin{equation}\label{Eq:HInfNorm}
\|S_{ii}\|_{\mathcal{H}_{\infty}} = \sup_{\omega \in \real} |S_{ii}(j\omega)| = |1 - \tfrac{1}{\beta}(\beta_i-b_i)|.
\end{equation}

%

\begin{figure}[ht!]
\begin{center}
\includegraphics[width=0.95\columnwidth]{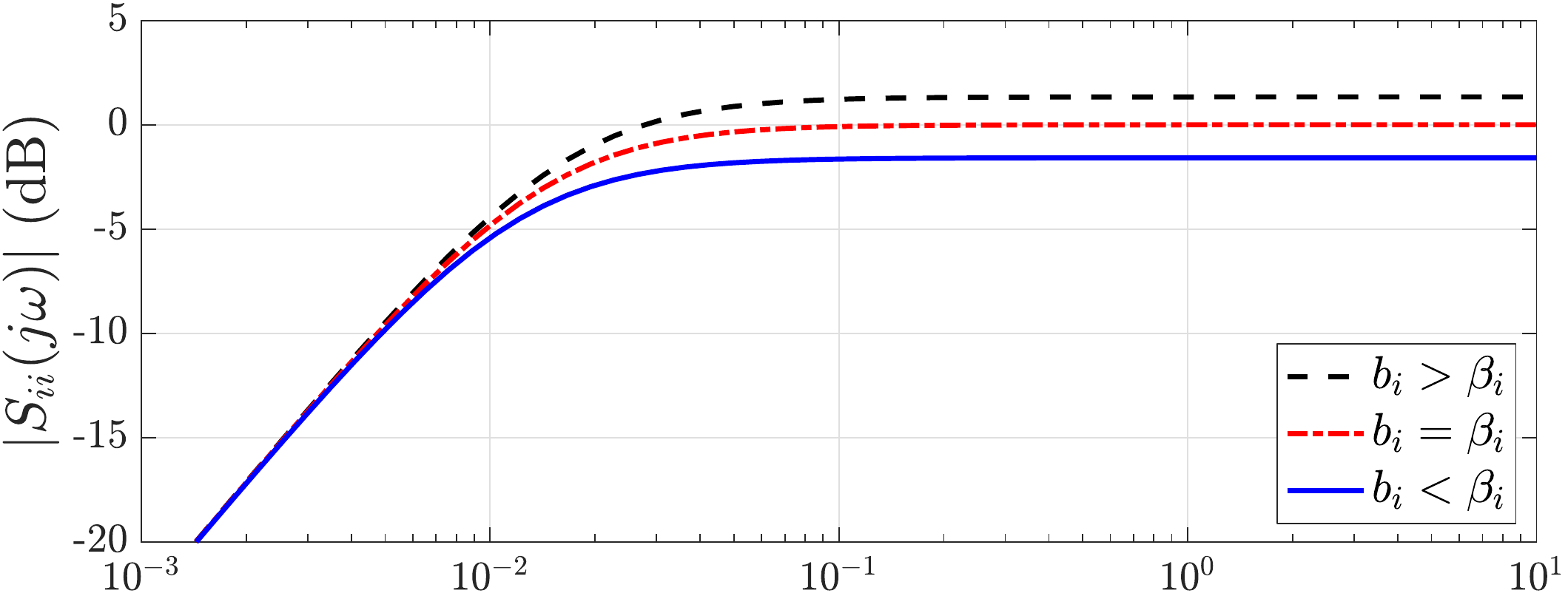}
\caption{Bode plot of $S_{ii}(s)$ for overbiased and underbiased tunings.}
\label{Fig:Bode}
\end{center}
\end{figure}

From either \eqref{Eq:HInfNorm} or from Figure \ref{Fig:Bode}, we see that overbiasing in area $k$ tends to \emph{increase} the peak sensitivity of the ACE in area $k$ with respect to local disturbances in area $k$; underbiasing has the opposite effect. Combining this information with our study from Section \ref{Sec:Overbiasing}, we conclude that there is a familiar tension in the tuning of AGC systems: an underbiased tuning improves local sensitivity to disturbances, but decreases the stability margin/convergence rate of the interconnected system. To the best of our knowledge, this is a new insight into the dynamic consequences of AGC bias tunings.
}

\section{Conclusions}
\label{Sec:Conclusion}

{We have derived a simple necessary and sufficient condition for diagonal stability of a class of matrices with a rank-1 interconnection structure. We extended the result to provide sufficient conditions for diagonal stability under higher-rank perturbations. Finally, we leveraged the main result to provide the first rigorous stability analysis of automatic generation control in interconnected power systems, {\tb and we examined some of the immediate tuning consequences of our analysis.} One open question is whether the higher-rank conditions in  \eqref{Eq:SigmaBound} can be refined,  for example, by recursively constructing diagonal Lyapunov functions as additional terms are considered in \eqref{Eq:GammaSVD}. To further develop the practical utility of our AGC analysis, future work will focus on extending the reduced dynamics \eqref{Eq:SlowTimeScaleAGC2} obtained for AGC to include turbine-governor nonlinearity and sampled-data implementation effects. 
}


\appendices
\section{Technical Lemmas and Proofs}
\label{App:1}

\begin{lemma}\label{Lem:1}
Let $A = \left[\begin{smallmatrix}A_{11} & A_{12}\\ 0 & A_{22}\end{smallmatrix}\right]$ be a block triangular matrix. Then $A$ is diagonally stable if and only if both $A_{11}$ and $A_{22}$ are diagonally stable.
\end{lemma}

\begin{pfof}{Lemma \ref{Lem:1}}
The ``only if'' direction is immediate. For the other direction, let $D_1, D_2 \succ 0$ be such that
\[
\begin{aligned}
A_{11}^{\T}D_1 + D_1A_{11} &\prec 0\\
A_{22}^{\T}D_2 + D_2 A_{22} &\prec 0.
\end{aligned}
\]
For $\epsilon > 0$ define $D_{\epsilon} = \mathrm{diag}(D_1,\epsilon D_2)$. Then
\[
A^{\T}D_{\epsilon} + D_{\epsilon}A = 
\begin{bmatrix}\epsilon(A_{11}^{\T}D_1 + D_1A_{11}) & \epsilon M\\
\star & A_{22}^{\T}D_2 + D_2A_{22}\end{bmatrix}
\]
where $M = D_{1}A_{12} + A_{12}^{\T}D_2$. As the (1,1) block is negative definite, the overall matrix is negative definite if and only if
\[
Q_2 - \epsilon M^{\T}(A_{11}^{\T}D_1 + D_1A_{11})^{-1}M \prec 0
\]
where $Q_2 =  A_{22}^{\T}D_2 + D_2A_{22} \prec 0$. This holds for sufficiently small $\epsilon$, which completes the proof.
\end{pfof}




\begin{lemma}\label{Lem:properties}
Let $A, P \in \real^{n \times n}$ with $P \succ 0$. The following statements are equivalent:
\begin{enumerate}
    \item 
    $A$ is diagonally stable with certificate $P \succ 0$;
    \item for any diagonal $\Delta_1, \Delta_2 \succ 0$ $\Delta_1 A \Delta_2$ is diagonally stable with certificate $\Delta_2P\Delta_1^{-1} \succ 0$;
    \item for any permutation matrix $\Pi$, $\Pi A \Pi^{\T}$ is diagonally stable with certificate $\Pi P \Pi^{\T}\succ 0$.
\end{enumerate}
\end{lemma}
\begin{pfof}{Lemma \ref{Lem:properties}}
The results easily follow from the definition of diagonal stability.
\end{pfof}



\bibliographystyle{IEEEtran}
\bibliography{brevalias,Main,JWSP,New}

\begin{thebibliography}{10}
\providecommand{\url}[1]{#1}
\csname url@samestyle\endcsname
\providecommand{\newblock}{\relax}
\providecommand{\bibinfo}[2]{#2}
\providecommand{\BIBentrySTDinterwordspacing}{\spaceskip=0pt\relax}
\providecommand{\BIBentryALTinterwordstretchfactor}{4}
\providecommand{\BIBentryALTinterwordspacing}{\spaceskip=\fontdimen2\font plus
\BIBentryALTinterwordstretchfactor\fontdimen3\font minus
  \fontdimen4\font\relax}
\providecommand{\BIBforeignlanguage}[2]{{%
\expandafter\ifx\csname l@#1\endcsname\relax
\typeout{** WARNING: IEEEtran.bst: No hyphenation pattern has been}%
\typeout{** loaded for the language `#1'. Using the pattern for}%
\typeout{** the default language instead.}%
\else
\language=\csname l@#1\endcsname
\fi
#2}}
\providecommand{\BIBdecl}{\relax}
\BIBdecl

\bibitem{RAH-CRJ:94}
R.~A. Horn and C.~R. Johnson, \emph{Topics in Matrix Analysis}.\hskip 1em plus
  0.5em minus 0.4em\relax Cambridge University Press, 1994.

\bibitem{EK-AB:00}
E.~Kaszkurewicz and A.~Bhaya, \emph{Matrix Diagonal Stability in Systems and
  Computation}.\hskip 1em plus 0.5em minus 0.4em\relax Springer, 2000.

\bibitem{LS-MA-EDS:10}
L.~{Scardovi}, M.~{Arcak}, and E.~D. {Sontag}, ``Synchronization of
  interconnected systems with applications to biochemical networks: An
  input-output approach,'' \emph{IEEE Trans. Autom. Control}, vol.~55, no.~6,
  pp. 1367--1379, 2010.

\bibitem{HKK:87}
H.~K. {Khalil}, ``Stability analysis of nonlinear multiparameter singularly
  perturbed systems,'' \emph{IEEE Trans. Autom. Control}, vol.~32, no.~3, pp.
  260--263, 1987.

\bibitem{RS-KSN:09}
R.~Shorten and K.~S. Narendra, ``On a theorem of redheffer concerning diagonal
  stability,'' \emph{Linear Algebra and its Applications}, vol. 431, no.~12,
  pp. 2317 -- 2329, 2009.

\bibitem{AR:15}
A.~Rantzer, ``Scalable control of positive systems,'' \emph{European Journal of
  Control}, vol.~24, pp. 72 -- 80, 2015.

\bibitem{MH:78}
P.~Moylan and D.~Hill, ``Stability criteria for large--scale systems,''
  \emph{IEEE Trans. Autom. Control}, vol.~23, no.~2, pp. 143--149, 1978.

\bibitem{DDS:78}
D.~D. {\v{S}}iljak, \emph{Large-scale systems: Stability and structure}.\hskip
  1em plus 0.5em minus 0.4em\relax North-Holland, 1978.

\bibitem{MV:81}
M.~Vidyasagar, \emph{Input-Output Analysis of Large-Scale Interconnected
  Systems: {D}ecomposition, Well-Posedness and Stability}.\hskip 1em plus 0.5em
  minus 0.4em\relax Springer Verlag, 1981.

\bibitem{MA-CM-AP:16}
M.~Arcak, C.~Meissen, and A.~Packard, \emph{Networks of Dissipative Systems:
  Compositional Certification of Stability, Performance, and Safety}.\hskip 1em
  plus 0.5em minus 0.4em\relax Springer Briefs in Control, Automation and
  Robotics, 2016.

\bibitem{MArcak:11}
M.~{Arcak}, ``Diagonal stability on cactus graphs and application to network
  stability analysis,'' \emph{IEEE Trans. Autom. Control}, vol.~56, no.~12, pp.
  2766--2777, 2011.

\bibitem{FD-SG:17}
F.~D{\"o}rfler and S.~Grammatico, ``Gather-and-broadcast frequency control in
  power systems,'' \emph{Automatica}, vol.~79, pp. 296 -- 305, 2017.

\bibitem{NM-CDP-AJvdS-JMAS:18}
N.~Monshizadeh, C.~De~Persis, A.~J. van~der Schaft, and J.~M.~A. Scherpen, ``A
  novel reduced model for electrical networks with constant power loads,''
  \emph{IEEE Trans. Autom. Control}, vol.~63, no.~5, pp. 1288--1299, 2018.

\bibitem{CDP-SG:20}
C.~De~Persis and S.~Grammatico, ``Continuous-time integral dynamics for a class
  of aggregative games with coupling constraints,'' \emph{IEEE Trans. Autom.
  Control}, vol.~65, no.~5, pp. 2171--2176, 2020.

\bibitem{XF-TA-MA-JTW-TB:06}
X.~Fan, T.~Alpcan, M.~Arcak, T.~Wen, and T.~Basar, ``A passivity approach
  to game-theoretic cdma power control,'' \emph{Automatica}, vol.~42, no.~11,
  pp. 1837--1847, 2006.

\bibitem{JTW-MA:04}
J.~T. Wen and M.~Arcak, ``A unifying passivity framework for network flow
  control,'' \emph{IEEE Trans. Autom. Control}, vol.~49, no.~2, pp. 162--174,
  2004.

\bibitem{HKK:02}
H.~K. Khalil, \emph{Nonlinear Systems}, 3rd~ed.\hskip 1em plus 0.5em minus
  0.4em\relax Prentice Hall, 2002.

\bibitem{AJW-BFW:96}
A.~J. Wood and B.~F. Wollenberg, \emph{Power Generation, Operation, and
  Control}, 2nd~ed.\hskip 1em plus 0.5em minus 0.4em\relax John Wiley \& Sons,
  1996.

\bibitem{PK:94}
P.~Kundur, \emph{Power System Stability and Control}.\hskip 1em plus 0.5em
  minus 0.4em\relax McGraw-Hill, 1994.

\bibitem{MDI-SL:96}
M.~D. Ili\'{c} and S.~Liu, \emph{Hierarchical Power Systems Control: Its Value
  in a Changing Industry}, ser. Advances in Industrial Control.\hskip 1em plus
  0.5em minus 0.4em\relax Springer, 1996.

\bibitem{MI-JZ:00}
M.~Ili\'{c} and J.~Zaborszky, \emph{Dynamics and Control of Large Electric
  Power Systems}.\hskip 1em plus 0.5em minus 0.4em\relax Wiley-IEEE Press,
  2000.

\bibitem{NERC:10bal}
NERC, ``Bal-005-1 balancing authority control,'' North American Electric
  Reliability Corporation, Tech. Rep., 2010.

\bibitem{RAH-CRJ:12}
R.~A. Horn and C.~R. Johnson, \emph{Matrix Analysis}, 2nd~ed.\hskip 1em plus
  0.5em minus 0.4em\relax Cambridge University Press, 1985.

\bibitem{SHL:14a}
S.~H. Low, ``Convex relaxation of optimal power flow, part i: Formulations and
  equivalence,'' \emph{IEEE Trans. Control Net. Syst.}, vol.~1, no.~1, pp.
  15--27, 2014.

\bibitem{NC:83}
N.~Cohn, ``The evolution of real time control applications to power systems,''
  in \emph{IFAC Symposium on Real Time Digital Control Applications}, vol.~16,
  no.~1, Guadalajara, Mexico, Jan. 1983, pp. 1 -- 17.

\bibitem{FPD-RJM-WFB:73a}
F.~P. {deMello}, R.~J. {Mills}, and W.~F. {B'Rells}, ``{A}utomatic {G}eneration
  {C}ontrol {P}art {I} -- {P}rocess {M}odeling,'' \emph{IEEE Trans. Power
  Apparatus \& Syst.}, vol. PAS-92, no.~2, pp. 710--715, 1973.

\bibitem{FPD-RJM-WFB:73b}
------, ``{A}utomatic {G}eneration {C}ontrol {P}art {II} -- {D}igital {C}ontrol
  {T}echniques,'' \emph{IEEE Trans. Power Apparatus \& Syst.}, vol. PAS-92,
  no.~2, pp. 716--724, 1973.

\bibitem{HGK-KCK-AB:75}
H.~G. {Kwatny}, K.~C. {Kalnitsky}, and A.~{Bhatt}, ``An optimal tracking
  approach to load-frequency control,'' \emph{IEEE Trans. Power Apparatus \&
  Syst.}, vol.~94, no.~5, pp. 1635--1643, 1975.

\bibitem{WBG:78}
W.~B. Gish, ``Automatic generation control -- notes and observations,'' Bureau
  of Reclamation, Engineering and Research Center, Tech. Rep. REC-ERC-78-6,
  1978.

\bibitem{HG-JS:80}
H.~Glavitsch and J.~Stoffel, ``Automatic generation control,'' \emph{Int. J.
  Electrical Power \& Energy Syst.}, vol.~2, no.~1, pp. 21 -- 28, 1980.

\bibitem{JC:85}
J.~Carpentier, ``To be or not to be modern, that is the question for automatic
  generation control (point of view of a utility engineer),'' \emph{Int. J.
  Electrical Power \& Energy Syst.}, vol.~7, no.~2, pp. 81 -- 91, 1985.

\bibitem{MSC:86}
M.~S. {Calovic}, ``Recent developments in decentralized control of generation
  and power flows,'' in \emph{Proc. IEEE CDC}, Athens, Greece, Dec. 1986, pp.
  1192--1197.

\bibitem{NJ-LSV-DNE-LHF-AGH:92}
N.~Jaleeli, L.~S. VanSlyck, D.~N. Ewart, L.~H. Fink, and A.~G. Hoffmann,
  ``Understanding automatic generation control,'' \emph{IEEE Trans. Power
  Syst.}, vol.~7, no.~3, pp. 1106--1122, 1992.

\bibitem{IEEE-Report-AGC:79}
M.~D. Anderson, ``Current operating problems associated with automatic
  generation control,'' \emph{IEEE Trans. Power Apparatus \& Syst.}, vol.~98,
  no.~1, pp. 88--96, 1979.

\bibitem{PKI-DPK:05}
P.~K. Ibraheem and D.~P. Kothari, ``Recent philosophies of automatic generation
  control strategies in power systems,'' \emph{IEEE Trans. Power Syst.},
  vol.~20, no.~1, pp. 346--357, 2005.

\bibitem{HHA-MEHG-RZ-EHF-PS:18}
H.~H. Alhelou, M.-E. Hamedani-Golshan, R.~Zamani, E.~Heydarian-Forushani, and
  P.~Siano, ``Challenges and opportunities of load frequency control in
  conventional, modern and future smart power systems: A comprehensive
  review,'' \emph{Energies}, vol.~11, no.~10, p. 2497, 2018.

\bibitem{DKM-FD-HS-SHL-SC-RB-JL:17}
D.~K. Molzahn, F.~D\"{o}rfler, H.~Sandberg, S.~H. Low, S.~Chakrabarti,
  R.~Baldick, and J.~Lavaei, ``A survey of distributed optimization and control
  algorithms for electric power systems,'' \emph{IEEE Trans. Smart Grid},
  vol.~8, no.~6, pp. 2941--2962, 2017.

\bibitem{FD-SB-JWSP-SG:19a}
F.~D\"{o}rfler, S.~Bolognani, J.~W. Simpson-Porco, and S.~Grammatico,
  ``Distributed control and optimization for autonomous power grids,'' in
  \emph{Proc. {ECC}}, Naples, Italy, Jun. 2019, pp. 2436--2453.

\bibitem{ED-NP:78}
E.~Davison and N.~Tripathi, ``The optimal decentralized control of a large
  power system: Load and frequency control,'' \emph{IEEE Trans. Autom.
  Control}, vol.~23, no.~2, pp. 312--325, 1978.

\bibitem{ANV-IAH-JBR-SJW:08}
A.~N. Venkat, I.~A. Hiskens, J.~B. Rawlings, and S.~J. Wright, ``Distributed
  mpc strategies with application to power system automatic generation
  control,'' \emph{IEEE Trans. Control Syst. Tech.}, vol.~16, no.~6, pp.
  1192--1206, 2008.

\bibitem{AM-LI-KU:16}
A.~M. Ersdal, L.~Imsland, and K.~Uhlen, ``Model predictive load-frequency
  control,'' \emph{IEEE Trans. Power Syst.}, vol.~31, no.~1, pp. 777--785,
  2016.

\bibitem{PRBM-PT:17}
P.~R.~B. Monasterios and P.~Trodden, ``Low-complexity distributed predictive
  automatic generation control with guaranteed properties,'' \emph{IEEE Trans.
  Smart Grid}, vol.~8, no.~6, pp. 3045--3054, 2017.

\bibitem{QL-MDI:12}
Q.~{Liu} and M.~D. {Ili\'{c}}, ``Enhanced automatic generation control
  ({E}-{AGC}) for future electric energy systems,'' in \emph{Proc. IEEE PESGM},
  2012, pp. 1--8.

\bibitem{DA:14}
D.~Apostolopoulou, ``Enhanced automatic generation control with uncertainty,''
  Ph.D. dissertation, University of Illinois at Urbana-Champaign, 2014.

\bibitem{CL-JHC:17}
C.~Lackner and J.~H. Chow, ``Wide-area automatic generation control between
  control regions with high renewable penetration,'' in \emph{IREP Symposium
  Bulk Power System Dynamics and Control}, 2017, pp. 1--10.

\bibitem{NL-CZ-LC:16}
N.~Li, C.~Zhao, and L.~Chen, ``Connecting automatic generation control and
  economic dispatch from an optimization view,'' \emph{IEEE Trans. Control Net.
  Syst.}, vol.~3, no.~3, pp. 254--264, Sep. 2016.

\bibitem{CZ-EM-SHL-JB:18}
C.~Zhao, E.~Mallada, S.~H. Low, and J.~Bialek, ``Distributed plug-and-play
  optimal generator and load control for power system frequency regulation,''
  \emph{Int. J. Electrical Power \& Energy Syst.}, vol. 101, pp. 1 -- 12, 2018.

\bibitem{JWSP:20b}
J.~W. Simpson-Porco, ``On stability of distributed-averaging
  proportional-integral frequency control in power systems,'' \emph{IEEE
  Control Syst. Let.}, vol.~5, no.~2, pp. 677--682, 2020.

\bibitem{DA-PWS-ADD:16}
D.~{Apostolopoulou}, P.~W. {Sauer}, and A.~D. {Dom\'{i}nguez-Garc\'{i}a},
  ``Balancing authority area model and its application to the design of
  adaptive {AGC} systems,'' \emph{IEEE Trans. Power Syst.}, vol.~31, no.~5, pp.
  3756--3764, 2016.

\bibitem{AAT-FZ-LX:11}
A.~A. {Thatte}, F.~Zhang, and L.~{Xie}, ``Frequency aware economic dispatch,''
  in \emph{Proc. NAPS}, Boston, MA, USA, Aug. 2011, pp. 1--7.

\bibitem{GZ-JM-QW:19}
G.~{Zhang}, J.~{McCalley}, and Q.~{Wang}, ``An {AGC} dynamics-constrained
  economic dispatch model,'' \emph{IEEE Trans. Power Syst.}, vol.~34, no.~5,
  pp. 3931--3940, 2019.

\bibitem{EE-ZT-JWSP-EF-MP-HH:20k}
E.~Ekomwenrenren, Z.~Tang, J.~W. Simpson-Porco, E.~Farantatos, M.~Patel, and
  H.~Hooshyar, ``Hierarchical coordinated fast frequency control using
  inverter-based resources,'' \emph{IEEE Trans. Power Syst.}, 2020,
  early-Access on on IEEE Xplore April 28th, 2021.

\bibitem{PC-SD-YCC-MP:20}
P.~Chakraborty, S.~Dhople, Y.~C. Chen, and M.~Parvania, ``Dynamics-aware
  continuous-time economic dispatch and optimal automatic generation control,''
  in \emph{Proc. {ACC}}, Denver, CO, USA, Jun. 2020, to appear.

\bibitem{CWT-RLC:76}
C.~W. {Taylor} and R.~L. {Cresap}, ``Real-time power system simulation for
  automatic generation control,'' \emph{IEEE Trans. Power Apparatus \& Syst.},
  vol.~95, no.~1, pp. 375--384, 1976.

\bibitem{JWSP:17b}
J.~W. Simpson-Porco, ``A theory of solvability for lossless power flow
  equations -- {P}art {II}: Conditions for radial networks,'' \emph{IEEE Trans.
  Control Net. Syst.}, vol.~5, no.~3, pp. 1373--1385, 2018.

\bibitem{fd-fb:09z}
F.~Dorfler and F.~Bullo, ``Synchronization and transient stability in power
  networks and non-uniform {K}uramoto oscillators,'' \emph{SIAM J Ctrl Optm},
  vol.~50, no.~3, pp. 1616--1642, 2012.

\bibitem{FD-FB:11d}
F.~D{\"o}rfler and F.~Bullo, ``Kron reduction of graphs with applications to
  electrical networks,'' \emph{IEEE Trans. Circuits \& Syst.~I}, vol.~60,
  no.~1, pp. 150--163, 2013.

\bibitem{ARB-DJH:81}
A.~R. Bergen and D.~J. Hill, ``A structure preserving model for power system
  stability analysis,'' \emph{IEEE Trans. Power Apparatus \& Syst.}, vol. 100,
  no.~1, pp. 25--35, 1981.

\bibitem{DJH-IMYM-90}
D.~J. Hill and I.~M.~Y. Mareels, ``Stability theory for differential/algebraic
  systems with application to power systems,'' \emph{IEEE Trans. Circuits \&
  Syst.}, vol.~37, no.~11, pp. 1416--1423, 1990.

\bibitem{AS-HKK:84}
A.~{Saberi} and H.~{Khalil}, ``Quadratic-type lyapunov functions for singularly
  perturbed systems,'' \emph{IEEE Trans. Autom. Control}, vol.~29, no.~6, pp.
  542--550, 1984.

\bibitem{NC:56}
N.~{Cohn}, ``Some aspects of tie-line bias control on interconnected power
  systems [includes discussion],'' \emph{Transactions of the American Institute
  of Electrical Engineers. Part III: Power Apparatus and Systems}, vol.~75,
  no.~3, pp. 1415--1436, 1956.

\bibitem{NERC:11}
N.~R. Subcommittee, ``Balancing and frequency control,'' North American
  Electric Reliability Corporation, Tech. Rep., 2011.

\bibitem{UCTE-App1:04}
UCTE, ``Ucte operating handbook appendix 1: Load-frequency control and area
  performance,'' 2004.

\end{thebibliography}

\end{document}